\documentclass[final,3p,times]{elsarticle}

\journal{Linear Algebra and its Applications}

\usepackage{amsthm}
\usepackage{array, makecell}
\usepackage{lipsum}
\usepackage{amsfonts}
\usepackage{graphicx}
\usepackage{epstopdf}
\usepackage{xcolor}
\usepackage{subcaption}
\usepackage{adjustbox}
\usepackage{amsmath}
\usepackage{algpseudocode}
\usepackage{algorithm}

\usepackage{xcolor}

\DeclareMathOperator*{\argmin}{arg\,min}

\usepackage{hyperref}
\usepackage{amsopn}
\DeclareMathOperator{\diag}{diag}

\begin{document}

\begin{frontmatter}
\title{Augmented Flexible Krylov Subspace methods with applications to Bayesian inverse problems}

\author{Malena Sabaté Landman\textsuperscript{1$\dagger$}\corref{CorrespondingAuthor}}
\ead{m.sabate.landman@emory.edu} 
\author{Jiahua Jiang \textsuperscript{2$\dagger$}}
\ead{j.jiang.3@bham.ac.uk}
\author{Jianru Zhang \textsuperscript{3}}
\author{Wuwei Ren \textsuperscript{3}}
\address{\textsuperscript{1$\dagger$} Department of Mathematics, Emory University, Atlanta, US}
\address{\textsuperscript{2$\dagger$}School of Mathematics, University of Birmingham, UK}
\address{\textsuperscript{3} School of Information Science and Technology, ShanghaiTech University}
\address{\textsuperscript{$\dagger$}These authors contributed equally to this work}

\vspace{10pt}

\begin{abstract}
This paper presents two new augmented flexible (AF)-Krylov subspace methods, AF-GMRES and AF-LSQR, to compute solutions of large-scale linear  discrete ill-posed problems that can be modeled as the sum of two independent random variables, exhibiting smooth and sparse stochastic characteristics respectively. Following a Bayesian modelling approach, this corresponds to adding a covariance-weighted quadratic term and a sparsity enforcing $\ell_1$ term in the original least-squares minimization scheme. To handle the $\ell_1$ regularization term, the proposed approach constructs a sequence approximating quadratic problems that are partially solved using augmented flexible Krylov–Tikhonov methods.

Compared to other traditional methods used to solve this minimization problem, such as those based on iteratively reweighted norm schemes, the new algorithms build a single (augmented, flexible) approximation (Krylov) subspace that encodes information about the different regularization terms through adaptable “preconditioning”. The solution space is then expanded as soon as a new problem within the sequence is defined. This also allows for the regularization parameters to be chosen on-the-fly at each iteration. Compared to most recent work on generalized flexible Krylov methods \cite{sdecompose}, our methods offer theoretical assurance of convergence and a more stable numerical performance. The efficiency of the new methods is shown through a variety of experiments, including a synthetic image deblurring problem, a synthetic atmospheric transport problem, and fluorescence molecular tomography reconstructions using both synthetic and real-world experimental data.
\end{abstract}


\cortext[CorrespondingAuthor]{Corresponding author}
\end{frontmatter}

\section{Introduction}
This paper is concerned with linear ill-posed inverse problems of the form
\begin{equation}\label{eqn:original}
 A u_{\text{true}} + e = b,\quad A \in \mathbb{R}^{m\times n},
\end{equation}
where $e \sim \mathcal{N}(0,R)$ is some unknown noise with (symmetric positive definite) covariance $R$ that affects the measurements $b$, and $u_{\text{true}}$ is the unknown object we want to reconstruct. Moreover, $A$ is assumed to be  ill-conditioned with ill-determined rank, i.e. the singular values of $A$ decay and cluster at zero without an evident gap between two consecutive ones. In this case, the choice of appropriate regularization is crucial to obtain a meaningful approximation of the solution, see e.g. \cite{Hanke1993RegularizationMF}. In particular, we are interested in solutions that can be modeled as the sum of two independent random variables with different priors:
\begin{equation} \label{eq:two_comp}
u_{\text{true}} = x + \xi.
\end{equation}
Following from the Bayesian interpretation of variational regularization as the maximum a posteriori (MAP) for the solution $u_{\text{true}}$, an appropriate regularization choice when considering solutions of the form \eqref{eq:two_comp} is
\begin{equation}\label{eqn:split_ls}
\min_{x, \xi} \left\{\| A (x + \xi) - b\|^2_{R^{-1}} + \lambda_x^2 R_x(x) + \lambda_\xi^2 R_{\xi}(\xi)\right\},
\end{equation}
where the regularization functionals $\lambda_x^2 R_x$ and $\lambda_\xi^2 R_\xi$ are related to the priors on $x$ and $\xi$, and where $\lambda_x, \lambda_\xi$ can also be interpreted as regularization parameters that must be determined. 
In this paper we focus on the case where $x$ is smooth and $\xi$ is sparse \cite{sdecompose}. That is, $x$ is modeled as a Gaussian random variable with mean $\mu_x \in \mathbb{R}^n$, covariance matrix $Q \in \mathbb{R}^{n\times n}$ (symmetric positive definite) and scaling parameter $\lambda_x$, and $\xi$ is modeled as a univariate Laplace distributed variable with mean  $[\mu_{\xi}]_j$ and scaling parameter $2{\lambda}_{\xi}^{-2}$, in the sense that
\begin{equation}\label{eq:priors}
    x \sim \mathcal{N}(\mu_x, {\lambda}_{x}^{-2} Q), \quad [\xi]_j \sim \mathcal{L}([\mu_{\xi}]_j, 2{\lambda}_{\xi}^{-2}), \quad 1 \leq j \leq n.
\end{equation}
Such scenarios arise in many imaging applications. For instance, in atmospheric inverse modeling,  anomalies correspond to sparse abnormally large values, while the average background exhibits a smooth behaviour. Similarly, in biomedical imaging problems like fluorescence molecular tomography, the reconstruction consists of sparsely distributed high-intensity fluorescent signals (typically concentrated in a tumor site) against a smooth background (e.g. healthy biological tissue). Given the priors in \eqref{eq:priors}, the MAP estimate of $u_{\text{true}}$ (given in terms of $x$ and $\xi$) is the solution of the following minimization problem:
\begin{equation}\label{eqn:problem}
\min_{x, \xi} \left\{ \| A (x + \xi) - b\|^2_{R^{-1}} + \lambda_x^2 \|x\|^2_{Q^{-1}} + \lambda_\xi^2 \|\xi\|_1 \right\}.
\end{equation}

For large-scale problems where $A$ is only accessible in the form of a function that efficiently computes matrix vector products between $A$ and, possibly, $A^{\top}$, on any given vector (i.e. $A$ does not have an exploitable structure and/or cannot be stored), the only way to solve problem \eqref{eqn:problem} is to apply an iterative method. Note that we assume that we can store a (small) number of basis vectors. 

The challenge of solving \eqref{eqn:problem} is two-fold. First, handling $\ell_1$ regularization poses computational challenges because of its non-differentiability at the origin.  Typically,  nonlinear optimization techniques or quadratic approximations within an inner-outer structure \cite{https://doi.org/10.1002/cpa.20303,IRN2,Arridge_2014,doi:10.1137/15M1037925} are employed to address this. Yet these approaches are very computationally expensive and can lead to slow convergence.  Furthermore, several other accelerated methods, including the split Bregman method \cite{goldstein2009split}, separable approximations \cite{3883}, and accelerations of the iterative shrinkage thresholding algorithm \cite{doi:10.1137/080716542}, require the upfront selection of various parameters, which can be a challenging task. A powerful alternative is the use of flexible Krylov methods
\cite{msl1,JulianneSilvia}, which can bypass the inner-outer scheme and also enable the automatic selection of the regularization parameters. Second, the covariance matrix $Q$ is usually very large and dense (e.g., where $Q$ is derived from a Mat$\Acute{\text{e}}$rn kernel or a dictionary collection), thus working with its inverse or square root can be cumbersome. To tackle this, the generalized Golub-Kahan process \cite{doi:10.1137/16M1081968} was introduced for problems with $\|\cdot\|_{Q^{-1}}$ regularization using a change of variables and reformulating the minimization problem. More recently, a solution decomposition hybrid projection method has been developed in \cite{sdecompose} to find approximate solutions of (\ref{eqn:problem}). The primary advantage of this approach is that the two different random variables can be reconstructed separately from the measurement. However, this method lacks convergence guarantees.  

In this paper, we propose a new augmented flexible Krylov subspace method to find the MAP estimate of the solution of Bayesian inverse problems in which the solution is presumed to be a combination of two random variables, one with smooth attributes and the other with sparse stochastic features. Our approach has four main components:  
\begin{enumerate}
\item Augmented and flexible versions of LSQR and GMRES, AF-LSQR and AF-GMRES, are presented. This is the first time that a GMRES-based solver is developed to solve regularization minimization problems involving both a weighted $\ell_2$ term and an $\ell_1$ term.

\item AF-LSQR and AF-GMRES are based on partial decompositions of $A$ (and possibly $A^{\top}$), providing a competitive and novel solution space that can successfully characterize the smooth and sparse features of the solution. 
\item By choosing appropriate approximations of \eqref{eqn:problem} that are (partially) solved at each iteration, theoretical convergence is guaranteed given fixed regularization parameters $\lambda_x$ and $\lambda_{\xi}$.
\item Different techniques to automatically and efficiently select the regularization parameters at each iteration for AF-GMRES and AF-LSQR are provided, with the understanding that other parameter choice criteria can also be seamlessly applied within this framework.
\end{enumerate}
The cornerstone of our methods is projecting a sequence of quadratic tangent majorants of (a smoothed version of) the original problem into a suitable space of increasing dimensions that can effectively represent solutions with both smooth and sparse elements. We demonstrate the efficacy of our algorithms across a range of large-scale Bayesian inverse problems, where the target solution contains both smooth and sparse stochastic features. In particular, in the last example with real-world experimental data from fluorescence molecular tomography, our AF-LSQR method yields satisfactory and stable results in a mere $40$ iterations. 

The paper is organised as follows. Section \ref{Chapter:method} presents a short review on iteratively reweighted methods, and introduces the new AF-LSQR and AF-GMRES methods. Heuristic techniques for determining regularization parameters are also provided in this section. Section \ref{Chapter:proof} describes the theoretical guarantees for AF-LSQR and AF-GMRES along with a succinct discussion on their interpretation. The results of various experiments are shown in Section \ref{Chapter:experiments}, and conclusions are provided in Section \ref{Chapter:conclusions}. 

\section{Augmented iteratively reweighted flexible Krylov subspace methods.}\label{Chapter:method}
This section presents two new iterative algorithms, AF-GMRES and AF-LSQR, to find approximated solutions to \eqref{eqn:problem}. First note that, for many problems, explicitly computing the inverse (or the square root) of $Q$ might be unfeasible. To alleviate this, and following \cite{sdecompose, doi:10.1137/16M1081968},  problem \eqref{eqn:problem} can be re-formulated as
\begin{equation}\label{eqn:def_phi_0}
\min_{x,\xi} \left\{\|A Q x + A \xi -b\|^2_{R^{-1}}+\lambda_{x}^2\|x\|^2_{Q}+\lambda_{\xi}^2 \| \xi\|_1\right\},
\end{equation}
where $Q$ is just required to be accessible in the form of matrix-vector products with any given vector. Second, an iteratively reweighted scheme is used to handle the minimization of the $\ell_1$ regularization term. This is explained in detail in the following subsection.

\subsection{Iteratively reweighted schemes}
Iteratively reweighted schemes involving an $\ell_1$-norm regularization term rely on the interpretation of the $\ell_1$-norm as a weighted $\ell_2$-norm, where the weights depend non-linearly on the solution. In particular, problem \eqref{eqn:def_phi_0} can be read as 
\begin{equation}\label{eq:reweighted}
\min_{x, \xi} \{\|AQx+ A\xi-b\|^2_{R^{-1}}+\lambda_x^2\|x\|_Q^2+\lambda_\xi^2\|W(\xi) \xi\|_2^2\},
\end{equation}
where $W(\xi)$ is a diagonal weighting that depends on $\xi$ as defined in \eqref{eq:weights}. To overcome the lack of smoothness of the functional in \eqref{eq:reweighted} at any vector with 0-valued components (and, in the discrete setting, to avoid numerical divisions by zero in that case), a more stable approximation of $W(\xi)$, named $\widetilde{W}^{(\tau)}(\xi)$, is used instead in practice:
\begin{equation}\label{eq:weights}
W(\xi) =\diag \left( (([\xi]_i)^{-\frac{1}{2}})_{i=1,...,n}\right)\approx \widetilde{W}^{(\tau)}(\xi)=\diag \left( (([\xi]^2_i+\tau^2)^{-\frac{1}{4}})_{i=1,...,n}\right),
\end{equation}
where $\tau$ is a parameter chosen ahead of the iterations. Note that using the weights $\widetilde{W}^{(\tau)}(\xi)$ instead of $W(\xi)$ in \eqref{eq:reweighted} leads to a smoothed version of problem \eqref{eq:reweighted}:
\begin{equation}\label{eqn:def_phi}
\min_{x,\xi} \phi (x,\xi) = \min_{x, \xi} \{\|AQx+ A\xi-b\|^2_{R^{-1}}+\lambda_x^2\|x\|_Q^2+\lambda_\xi^2\|\widetilde{W}^{(\tau)}(\xi) \xi\|_2^2\}.
\end{equation}
Formally, the original problem \eqref{eq:reweighted} can be recovered from \eqref{eqn:def_phi} setting $\tau = 0$. A well-established framework to solve problem \eqref{eqn:def_phi} is the use of a particular family of majorization-minimization schemes, also known as iteratively reweighted least squares (IRLS) \cite{https://doi.org/10.1002/cpa.20303} or iteratively reweighted norm (IRN) \cite{IRN2}. They consist on the local approximation of $\phi (x,\xi)$ by a sequence of quadratic functionals $\phi_k (x,\xi)$ so that \eqref{eqn:def_phi} can be approximated by a sequence of problems of the form
\begin{equation}\label{eqn:def_phi_k}
\min_{x,\xi} \phi_k(x,\xi)=\min_{x, \xi} \{\|AQx+ A\xi-b\|^2_{R^{-1}}+\lambda_{x}^2\|x\|_Q^2+\lambda_{\xi}^2\|W_k \xi\|_2^2 + c_k\},
\end{equation}
where $W_k = \widetilde{W}^{(\tau)}(\xi_{k-1})$, $c_k$ is an iteration-dependent constant term with respect to $x$ and $\xi$ and $\lambda_x$ and $\lambda_\xi$ have absorbed other possible multiplicative constants. Note that $c_k$ has to be chosen appropriately such that $\phi_k(x, \xi)$ is a quadratic tangent majorant of $\phi (x,\xi)$ at $\xi = \xi_{k-1}$. By definition, this implies that $\phi_k(x, \xi) \geq \phi(x, \xi)$ for all $x, \xi \in \mathbb{R}^n$, and that $\phi_k(x, \xi_{k-1}) = \phi(x, \xi_{k-1})$ and $\nabla \phi_k (x,\xi_{k-1}) = \nabla \phi(x, \xi_{k-1})$ for all $x \in \mathbb{R}^n$, see e.g. \cite{Huang2017, msl1}. Moreover, problems of the form \eqref{eqn:def_phi_k} can be transformed into standard form, so that the weights $W^{-1}_k$ can be interpreted as solution-dependent right preconditioners for the sparse component of the solution of:
\begin{equation}\label{eq:std_form}
\min_{y, \xi} \{\|AQx+AW^{-1}_k \xi-b\|^2_{R^{-1}}+\lambda_x^2\|x\|_Q^2+\lambda_{\xi}^2\|\xi\|_2^2\},
\end{equation}
providing a natural setting for flexible Krylov methods. Note that, as mentioned in  \cite{msl1}, the use of the preconditioning inspired by expression \eqref{eq:std_form} generates a solution space that incorporates prior information and adapts to the solution throughout the iterations. However, flexible Krylov methods do not solve a projection of the original problem in general and, differently from explicitly iteratively reweighted methods, might not converge to the solution of the smoothed version of the original problem \eqref{eqn:def_phi_0} defined in \eqref{eqn:def_phi}. For this reason, reweighting the regularization term concerning $\xi$ as in \eqref{eqn:def_phi_k} at every iteration is crucial to prove convergence of the solution computed using augmented flexible Krylov methods to the solution of the smoothed version of the original problem as will be explained in Section \ref{Chapter:proof}.

\subsection{Augmented flexible Arnoldi decomposition and AF-GMRES}
In this section a modified Arnoldi decomposition is proposed that includes both basis vectors with flexible preconditioning and basis vectors with fixed preconditioning. Given a square system matrix $A \in \mathbb{R}^{n \times n}$, right iteration-dependent preconditioning matrices $W^{-1}_k \in \mathbb{R}^{n \times n}$ for $\xi$, and right iteration-independent preconditioning $Q\in \mathbb{R}^{n \times n}$ for $x$, the following expression involving $A$ is updated at each iteration $k$ with two new columns of $\hat{Z}_{2k-1}$, $\hat{V}_{2k}$ and $\hat{H}_{2k,2k-1}$:
\small

\begin{equation}\label{eqn:AFArnoldi}
A 
\underbrace{
[Q v_1 \,\, z_1 \,\, Q v_2 \,\,z_2 \,\,Q v_3, ...]
}_{\hat{Z}_{2k-1}}=
\underbrace{
[v_1  \,\, v_2 \,\, \bar{v}_2 \,\,v_3 \,\,\bar{v}_3\,\,  v_4
, ... ]}_{\hat{V}_{2k}}
\underbrace{
\begin{bmatrix}
h_{11}   & h_{12}  & h_{13}   & h_{14}  & h_{15}  & ... \\
h_{21} & h_{22}  & h_{23}   & h_{24}  & h_{25} & ...  \\
           & h_{32} & h_{33}   & h_{34}  & h_{35}   & ... \\
           &           & h_{43}& h_{44}  & h_{45}  & ...  \\
           &           &            & h_{54}  & h_{55}  & ...    \\
           &           &            &           & h_{65}  & ... \\
           &           &            &           &     & ...        
\end{bmatrix}}_{\hat{H}_{2k,2k-1}},
\end{equation}
\normalsize
where $v_1=b/\|b\|_{R^{-1}}$ and $z_i= W_i^{-1} v_i$. Notation-wise, a single subscript in a matrix indicates the number of columns instead of the number of iterations, while two subscripts indicate the matrix dimensions. Here, $\hat{V}_{2k} \in \mathbb{R}^{n \times {2k}}$ has $R^{-1}$-orthogonal columns by construction: 
\begin{equation}\label{eq_Rorth}
    \hat{V}^{\top}_{2k} R^{-1} \hat{V}_{2k}=I_{2k}.
\end{equation}
Note that one has to be careful in handling potential breakdown due to linear dependence of subspace vectors. This can be avoided numerically checking for breakdown and avoiding adding that vector in the solution subspace. Since two vectors are added to the solution subspace at each iteration, this does not necessarily imply algorithmic break-down. 

Finally, note that we can define $ H^{(1)} \in  \mathbb{R}^{2k \times k}$ (resp. $ H^{(2)} \in  \mathbb{R}^{2k \times k-1}$) as the sub-matrices containing only the odd (resp. even) columns of the matrix $ \hat{H}_{2k,2k-1}$, so that:
\begin{equation}\label{eqn:AFArnoldi_split}
A Q V_k = \hat{V}_{2k} H^{(1)}, \quad A Z_{k-1} = \hat{V}_{2k} H^{(2)}.
\end{equation}
Here $Z_{k-1} \in  \mathbb{R}^{n \times {2k-1}}$ has one column less than $V_k$ 
due to the initialization of the spaces and the fact that we want to update the flexible preconditioner as soon as a new solution is available.

Augmented flexible (AF-) GMRES can be proposed as a generalization of GMRES using expressions \eqref{eqn:AFArnoldi} and \eqref{eqn:AFArnoldi_split} to provide a projection of \eqref{eqn:def_phi_k} onto a space of increasing dimensions. The solution is then computed at each iteration  $k$, for 1 $\leq k \leq n/2$, as:
\begin{equation}\label{eq:GMRES_sol1}
u_k= x_k + \xi_k = \hat{Z}_{2k-1}y_{2k-1}=Q V_{k}y_{2k-1}^{(1)}+Z_{k-1}y_{2k-1}^{(2)}
\end{equation}
where 
\begin{equation}\label{eqn:projected_AFGMRES_1}
y_{2k-1}=\argmin_{y} \left\{\| A \hat{Z}_{2k-1} \,
y
- b \|_{R^{-1}}^2  + \lambda_x^2 \| V_k y^{(1)}\|_Q^2+ \lambda_\xi^2 \|W_{k} Z_{k-1}y^{(2)}\|_2^2\right\}
\end{equation}
and $y^{(1)} \in \mathbb{R}^k$ (resp. $y^{(2)} \in \mathbb{R}^{k-1}$) correspond to the odd (resp. even) components of $y$. To project the second term we need to $Q-$reorthogonalise the columns of $V_k$, i.e., finding
\begin{equation}\label{eq:Qorth}
V_k = \widetilde{V}_k \widetilde{H}_{k,k} \quad \text{such that} \quad {\widetilde{V}}^{T}_{k} Q \widetilde{V}_k = I_{k,k}.
\end{equation}
Note that this can be updated throughout the iterations and it does not require extra matrix-vector products with $Q$ (lines 10-12 of Algorithm \ref{alg:AFGMRES}). Then, the first component of the solution in \eqref{eq:GMRES_sol1} can be expressed as
$V_k \, y^{(1)}= \widetilde{V}_k \, \widetilde{H}_{k,k}\, y^{(1)}$. Moreover, the QR factorization of the tall and skinny projected regularization matrix, i.e. 
\begin{equation}\label{eq_qr_WZ}
Q^{(wz)}R^{(wz)}=W_{k}Z_{k-1},   
\end{equation}
can be computed efficiently at each iteration, see e.g. \cite{doi:10.1137/080731992}. Therefore, $y_{2k-1}$ in \eqref{eq:GMRES_sol1} is computed by projecting \eqref{eqn:projected_AFGMRES_1} using \eqref{eq_Rorth},\eqref{eqn:AFArnoldi_split}, \eqref{eq:Qorth} and \eqref{eq_qr_WZ}, so that:
\begin{equation}\label{eqn:projected_AFGMRES_2}
y_{2k-1} = \argmin_{y} \left\{\| {H}_{2k+1,2k} \,
y
- \textcolor{black}{\|b\|_{R^{-1}}  e_1}\|_2^2  + \lambda_x^2 \| \widetilde{H}_{k,k} y^{(1)}\|_2^2+ \lambda_\xi^2 \|R^{(wz)} y^{(2)}\|_2^2\right\}.
\end{equation}
Assuming no algorithmic break-down has happened, the AF-GMRES method can be extended after $\hat{V}_{2k}$ has $n$ columns by considering $(x_k, \xi_k) = \argmin_{x,\xi\in \mathbb{R}^n} \phi_k(x,\xi)$, where $\phi_k(x,\xi)$ is defined in \eqref{eqn:def_phi_k}. Since $\phi_k(x,\xi)$ is a quadratic tangent majorant of $\phi(x,\xi)$ at $(x_{k-1},\xi_{k-1})$, this corresponds to an iteration of a majorization minimization (MM) scheme for minimizing $\phi(x,\xi)$ as defined in \eqref{eqn:def_phi}. This situation is not expected to happen in practice for large-scale problems; however, since $A$ is ill-posed one can expect most relevant directions in the range of $A$ to be represented in the Krylov subspace after a small number of iterations \cite{Krylov_props}. In other other words, after a small amount of iterations, AF-GMRES mimics an MM scheme.  
Solid theoretical foundations for the method can be found in Section \ref{Chapter:proof}. A detailed implementation of this is given in Algorithm \ref{alg:AFGMRES}.

\begin{algorithm}
\caption{Augmented Flexible (AF-) GMRES} \label{alg:AFGMRES}
\begin{algorithmic}[1]
\State Input: $A$, $Q$, $R$, $b$.
\State Initialize: 
\State \ \ \ \ \ Compute $\hat{v}_1 = v_1 = b /\|b\|_{R^{-1}}$
\State \ \ \ \ \ Compute $\tilde{v}_1 = \hat{v}_1 /\|\hat{v}_1\|_{Q}$, \quad $\tilde{h}_{1,1}=\|\hat{v}_1\|_{Q}$
\State  \ \ \ \ \  Compute $\hat{z}_1 = Q \, \hat{v}_1$ 
\For{$k=1,\dots$}
\State  $w = A \hat{z}_{2k-1}$ \Comment{mvp with $A$}
\State \textbf{for} $i$=1\,...\,2$k$-1 \textbf{do} \quad $h_{i,2k-1} =\hat{v_i}^{\top} (R^{-1}\, w)$, \, \, $w = w-h_{i,2k-1}\, \hat{v}_{i}$ \, \textbf{end for}
\State  Compute $\hat{v}_{2k} = v_{k+1} = w / \|w\|_{R^{-1}}$, \quad $h_{2k,2k-1}=\|w\|_{R^{-1}}$
\State $w = v_{k+1}$
\State \textbf{for} $i$=1\,...\,$k$ \textbf{do} \quad $\tilde{h}_{i,k+1} =v_i^{\top} (Q\, w)$, \, \, $w = w-\tilde{h}_{i,k+1}\, v_{i}$ \, \textbf{end for}
\State  Compute $\tilde{v}_{k+1} = w /\|w\|_{Q}$, \quad $\tilde{h}_{k+1,k+1} = \|w\|_{Q}$
\State Compute $y^{(1)}$ and $y^{(2)}$ in \eqref{projected_trick} using suitable regularization param. criteria 
\State Update $x_k$ and $\xi_k$ in \eqref{eq:GMRES_sol1} \Comment{update solution}
\State Update weights $W_k=\widetilde{W}^{(\tau)}(\xi_k)$ according to \eqref{eq:weights} after initialization
\State Compute $\hat{z}_{2k} = W_k^{-1} \hat{v}_{2k-1}$
\State $w=A \, \hat{z}_{2k} $ \Comment{ mvp with $A$}
\State \textbf{for} $i$=1\,...\,2$k$ \textbf{do} \quad $h_{i,2k} =\hat{v_i}^{\top} (R^{-1}\, w)$, \, \, $w = w-h_{i,2k}\, \hat{v}_{i}$ \textbf{end for}
\State  Compute $\hat{v}_{2k+1} =  w / \|w\|_{R^{-1}}$, \quad $h_{2k+1,2k}=\|w\|_{R^{-1}}$
\State  Compute $\hat{z}_{2k+1} =  Q \, \hat{v}_2k$
\EndFor
\end{algorithmic}
\end{algorithm}

\subsection{Augmented flexible Golub-Kahan decomposition and AF-LSQR}

Similarly to the previous section, a modified  Golub-Kahan decomposition is proposed, so that it includes both basis vectors with flexible preconditioning $W^{-1}_k$ for $\xi$ and basis vectors with fixed preconditioning $Q$ for $x$. At each iteration of the method, the following partial decompositions of $A$ and $A^{\top}$ are updated:

\begin{equation}\nonumber
A
\underbrace{
[Q v_1 \,\, z_1 \,\, Q v_2 \,\,z_2 \,\, Q v_3, ...  ]
}_{\hat{Z}_{2k-1}}=
\underbrace{
[u_1  \,\, u_2 \,\, \bar{u}_2 \,\,u_3 \,\,\bar{u}_3\,\,  u_4, ... 
]}_{\hat{U}_{2k}}
\underbrace{
\begin{bmatrix}
m_{11}   & m_{12}  & m_{13}   & m_{14}  & m_{15}   & ... \\
m_{21} & m_{22}  & m_{23}   & m_{24}  & m_{25} & ...  \\
           & m_{32} & m_{33}   & m_{34}  & m_{35}  & ...  \\
           &           & m_{43} & m_{44}  & m_{45} & ...   \\
           &           &            & m_{54} & m_{55}  & ...    \\
           &           &            &           & m_{65} & ...  \\
           &           &            &           &           & ...  
\end{bmatrix}}_{\hat{M}_{2k,2k-1}},
\end{equation}
\begin{equation}\label{GK2_2}
A^{\top} R^{-1}
\underbrace{
[u_1  \,\, u_2 \,\, \bar{u}_2 \,\,u_3 \,\,\bar{u}_3, ... ]
}_{\hat{U}_{2k-1}}=
\underbrace{
[v_1  \,\, v_2 \,\, \bar{v}_2 \,\,v_3 \,\,\bar{v}_3, ... ]}_{\hat{V}_{2k-1}}
\underbrace{
\begin{bmatrix}
t_{11}  & t_{12}  & t_{13}   & t_{14}  & t_{15} & ...   \\
 & t_{22}& t_{23}   & t_{24}  & t_{25} & ...  \\
           &  &  t_{33}  & t_{34}  & t_{35}  & ...  \\
           &           &    & t_{44}   & t_{45} & ...   \\
           &           &            &  & t_{55}   & ...    \\
           &           &            &           &         & ...    
\end{bmatrix}}_{\hat{T}_{2k-1,2k-1}}.
\end{equation}
Here, $R^{-1}$-orthogonality is imposed between the columns of $\hat{U}_{2k}$, and $Q$-orthogonality is imposed between the columns of $\hat{V}_{2k-1}$ by construction, i.e.
\begin{equation}\label{oth_prop}
\hat{U}^{\top}_{2k} R^{-1} \hat{U}_{2k}=I_{2k} \quad \text{and} \quad \hat{V}^{\top}_{2k-1} Q \hat{V}_{2k-1}=I_{2k-1}.
\end{equation}
Note that we can define $M^{(1)} \in  \mathbb{R}^{2k \times k}$ (resp. $ M^{(2)} \in  \mathbb{R}^{2k \times k-1}$) as the sub-matrices containing only the odd (resp. even) columns of the matrix $ \hat{M}_{2k,2k-1}$, so that:
\begin{equation}\label{eqn:AFGK_split}
A Q V_k = \hat{U}_{2k} M^{(1)}, \quad A Z_{k-1} = \hat{U}_{2k} M^{(2)}.
\end{equation}

The new AF-LSQR is defined so that at each iteration $1 \leq k \leq n/2$ an approximation $x_k \in \mathcal{R}(\hat{Z}_{2k-1})$ to the solution of \eqref{eqn:def_phi} is computed as:
\begin{equation}\label{eq:x_and_xi}
u_k = x_k + \xi_k =\hat{Z}_{2k-1}y_{2k-1}=Q V_{k}y_{2k-1}^{(1)}+Z_{k-1}y_{2k-1}^{(2)}
\end{equation}
where
\begin{equation}\label{eqn:AF-LSQR_2}
y_{2k-1}=\argmin_{y} \left\{\| {M}_{2k,2k-1} \,
y
- \|b\|_{R^{-1}}  e_1\|_2^2  + \lambda_x^2 \|y^{(1)}\|_2^2+ \lambda_{\xi}^2 \|W_{k} Z_{k-1}y^{(2)}\|_2^2\right\}.
\end{equation}
Here $y^{(1)} \in \mathbb{R}^k$ (resp. $y^{(2)} \in \mathbb{R}^{k-1}$) correspond to the odd (resp. even) components of $y$. 
Analogously to the case of AF-GMRES, the QR factorization of the regularization matrix, i.e. $Q^{(wz)}R^{(wz)}=W_{k}Z_{k-1}$, can be computed efficiently at each iteration, see e.g. \cite{doi:10.1137/080731992} so that \eqref{eqn:AF-LSQR_2} is equivalent to
\begin{equation}\label{eqn:projected_AFLSQR_1}
y_{2k-1}=\argmin_{y} \left\{ \| {M}_{2k,2k-1} \,
y
- \|b\|_{R^{-1}}  e_1\|_2^2  + \lambda_x^2 \|y^{(1)}\|_2^2+ \lambda_{\xi}^2 \|R^{(wz)} y^{(2)}\|_2^2 \right\}.
\end{equation}

Moreover, the AF-LSQR method can be extended after $\hat{U}_{2k}$ has $n$ columns by considering $(x_k, \xi_k)$=$\argmin_{x,\xi\in \mathbb{R}^n} \phi_k(x,\xi)$, where $\phi_k(x,\xi)$ is defined in \eqref{eqn:def_phi_k}. 
Solid theoretical foundations for the method can be found in Section \ref{Chapter:proof} and a detailed implementation of this is given in Algorithm \ref{alg:AFLSQR}.

\begin{algorithm}
\caption{Augmented Flexible (AF-) LSQR} \label{alg:AFLSQR}
\begin{algorithmic}[1]
\State Input: $A$, $Q$, $R$, $b$.
\State Initialize: 
\State \ \ \ \ \ Compute $\hat{u}_1 = b/ \|b\|_{R^{-1}}$
\State  \ \ \ \ \  $v=A^{\top} R^{-1} \, \hat{u}_1$, \quad $\hat{v}_1 =v / \|v\|_{Q}$, \quad $t_{1,1} =\|v\|_{Q}$, \quad $\hat{z}_1 =Q \hat{v}_1$
\State   \ \ \ \ \ $u = A \hat{z}_1$, $m_{1,1}=\hat{u}_1^{\top} (R^{-1}  u)$,  $u = u-m_{1,1}\hat{u}_1$, $\hat{v}_2 = u /\|u\|_{R^{-1}}$, $m_{2,1} = \|u\|_{R^{-1}}$ 
\State   \ \ \ \ \ Update solution $u_1$ according to \eqref{projected_trick} and \eqref{eq:x_and_xi}, with a suitable reg. param.
\Statex  \ \ \ \ \ criteria.  and initialize weights $W_1=\widetilde{W}^{(\tau)}(u_1)$ according to \eqref{eq:weights} 
\State  \ \ \ \ \  Compute $\hat{z}_2 = W_1^{-1} \, \hat{v}_1$, \quad $u = A \hat{z}_2$
\State  \ \ \ \ \ $m_{1,2}=\hat{u}_1^{\top} (R^{-1}\, u)$, \, $m_{2,2}=\hat{u}_2^{\top} (R^{-1}\, u)$, \, $u = u - m_{1,2}\hat{u}_1 - m_{2,2}\hat{u}_2$
\State  \ \ \ \ \ Compute $\hat{u}_2 = u / \|u\|_{R^{-1}}$, \quad $m_{3,2}= \|u\|_{R^{-1}}$
\For{$k=2,\dots$} 
\State $v=A^{\top} R^{-1} \, \hat{u}_{2k-2}$ \Comment{ mvp with $A^{\top}$}
\State \textbf{for} $i$=1\,...\,2$k$-3 \textbf{do} \quad $t_{i,2k-2} =\hat{v_i}^{\top} (Q\, v)$, \, \, $v = v-t_{i,2k-2}\, \hat{v}_{i}$ \textbf{end for}
\State  Compute $\hat{v}_{2k-2} =  v / \|v\|_Q$, \quad $t_{2k-2,2k-2}=\|v\|_Q$
\State  Compute $\hat{z}_{2k-1} = Q \, \hat{v}_{2k-2} $
\State $v=A^{\top} R^{-1} \, \hat{u}_{2k-1}$ \Comment{ mvp with $A^{\top}$}
\State \textbf{for} $i$=1\,...\,2$k$-2 \textbf{do} \quad $t_{i,2k-1} =\hat{v}_i^{\top} (Q\, v)$, \, \, $v = v-t_{i,2k-1}\, \hat{v}_{i}$ \textbf{end for}
\State  Compute $\hat{v}_{2k-1} =  v / \|v\|_Q$, \quad $t_{2k-1,2k-1}=\|v\|_Q$
\State $u=A \hat{z}_{2k-1}$ \Comment{ mvp with $A$}
\State \textbf{for} $i$=1\,...\,2$k$-1 \textbf{do} \quad $m_{i,2k-1} =\hat{u_i}^{\top} (R^{-1}\, u)$, \, \, $u = u-m_{i,2k-1}\, \hat{u}_{i}$ \textbf{end for}
\State  Compute $\hat{u}_{2k} =  u / \|u\|_{R^{-1}}$, \quad $m_{2k,2k-1}=\|u\|_{R^{-1}}$
\State Compute $y^{(1)}$ and $y^{(2)}$ in \eqref{projected_trick} using suitable regularization param. criteria 
\State Update $x_k$ and $\xi_k$ in \eqref{eq:x_and_xi} \Comment{update solution}
\State Update weights $W_k=\widetilde{W}^{(\tau)}(\xi_k)$ according to \eqref{eq:weights}
\State Compute $\hat{z}_{2k} = W_k^{-1} \hat{v}_{2k-2}$
\State $u=A \, \hat{z}_{2k} $ \Comment{ mvp with $A$}
\State \textbf{for} $i$=1\,...\,2$k$ \textbf{do} \quad $m_{i,2k} =\hat{u_i}^{\top} (R^{-1}\, u)$, \, \, $u = u-m_{i,2k}\, \hat{u}_{i}$ \textbf{end for}
\State  Compute $\hat{u}_{2k+1} =  u / \|u\|_{R^{-1}}$, \quad $m_{2k+1,2k}=\|u\|_{R^{-1}}$

\EndFor
\end{algorithmic}
\end{algorithm}

\subsection{Computational considerations}
Note that the projected problems \eqref{eqn:projected_AFGMRES_2} and \eqref{eqn:projected_AFLSQR_1} can be computed at each iteration using the following expression: 
\small
\begin{equation} \label{projected_trick}
\argmin_{y^{(1)},y^{(2)}}
\left\| 
\begin{bmatrix}
A V_k &  A Z_{k-1}  \\
\lambda_x L & 0_{k,k-1} \\
0_{k-1,k}  & \lambda_\xi R^{(wz)} 
\end{bmatrix}
\begin{bmatrix}
y^{(1)}\\
y^{(2)}
\end{bmatrix} -
\|b\|_{R^{-1}} e_1 
\right\|_2^2 =
\argmin_{y^{(1)},y^{(2)}}
\left\| 
\begin{bmatrix}
K^{(1)} & K^{(2)} \\
\lambda_x L & 0_{k,k-1} \\
0_{k-1,k}  & \lambda_\xi R^{(wz)} 
\end{bmatrix} 
\begin{bmatrix}
y^{(1)}\\
y^{(2)}
\end{bmatrix} -
\|b\|_{R^{-1}} e_1
\right\|_2^2,
\end{equation}
\normalsize
where $K^{(1)}$, $K^{(2)}$, $L$ correspond to $H^{(1)}$, $H^{(2)}$, $\widetilde{H}_{k,k}$ defined in \eqref{eqn:AFArnoldi_split} for AF-GMRES and to $M^{(1)}$, $M^{(2)}$, $I_{k,k}$ defined in \eqref{eqn:AFGK_split} for AF-LSQR. 

\subsection{Regularization parameter and initialization choices}
Choosing suitable regularization parameters ${\lambda}_x$ and ${\lambda}_\xi$ is crucial for finding a meaningful solution of (\ref{eqn:def_phi_k}). However, this can be a challenging task for large-scale problems. Alternatively, projections methods based on explicitly updating partial factorizations of the system matrix, such as AF-GMRES and the AF-LSQR, allow for the regularization parameters to be changed at each iteration without increasing the computational cost. This motivates the choice of regularization parameters that are suitable for each of the projected problems independently, so that, effectively, the regularization parameters in \eqref{projected_trick} have an explicit dependence on the iteration number $k$, namely ${\lambda}^{(k)}_x$ and ${\lambda}^{(k)}_\xi$. This is, more generally, also the framework of hybrid methods \cite{hybrid_review}, and suitable extensions to two regularization parameters are studied in \cite{sdecompose}. In this section we introduce three techniques to determine the values of ${\lambda}^{(k)}_x$ and ${\lambda}^{(k)}_\xi$ for the projected problems. Define $K = \begin{bmatrix}K^{(1)} & K^{(2)}\end{bmatrix}$. The minimizer of \eqref{projected_trick} is given by
\begin{equation}\label{eq:C}
    y_k (\lambda_x,\lambda_{\xi}) = \begin{bmatrix}
        y_k^{(1)} \\ y_k^{(2)}
    \end{bmatrix}= \underbrace{\left(K^{\top} K + \begin{bmatrix}
{\lambda}^{2}_{x} L^{\top} L & 0_{k,k-1} \\
0_{k-1,k}  & {\lambda}^{2}_{\xi} {R^{(wz)}}^{T} R^{(wz)} 
\end{bmatrix}\right)^{-1} K^{\top}}_{\large C({\lambda}_{x}, {\lambda}_{\xi})} \|b\|_{R^{-1}} e_1 .
\end{equation}

Just as a reference measure to evaluate the performance of the studied parameter selection methods, we define the optimal parameters as 
\begin{equation}
\label{eq:opt}
    (\lambda_x^{(k)},\lambda_{\xi}^{(k)} )=\argmin_{{\lambda}_{x}, {\lambda}_{\xi}} \| u_k(\lambda_x,\lambda_{\xi}) - u_{\text{true}} \|^2_2,
\end{equation}
where $u_{\text{true}}$ is the true solution that is not available in practice. As previously noted, this definition is primarily used to gauge the effectiveness of the parameter selection methods. If (a good estimation of) the noise level $\sigma$ is known a priori, the discrepancy principle (DP) can be used to select the regularization parameters for the projected problem, where 
\begin{equation}\label{eq:dp}
     (\lambda_x^{(k)},\lambda_{\xi}^{(k)}) =\argmin_{{\lambda}_{x}, {\lambda}_{\xi}}\left| \|K y_k({\lambda}_{x}, {\lambda}_{\xi}) -\|b\|_{R^{-1}} e_1\|_2^2 - \tau_{dp} n \sigma^2 \ \right| ,
\end{equation}
where $\tau_{dp} \geq 1$ is a constant close to 1, also referred to as the safety factor. If a-priori knowledge of the noise level is not available, the weighted generalized cross validation (WGCV) method can be adapted instead, where the regularization parameters $\lambda_x,\lambda_{\xi}$ are selected such that 
\begin{equation}\label{eq:wgcv}
    (\lambda_x^{(k)},\lambda_{\xi}^{(k)}) =\argmin_{{\lambda}_{x}, {\lambda}_{\xi}} \frac{\|K y_k({\lambda}_{x}, {\lambda}_{\xi}) - \|b\|_{R^{-1}} e_1\|_2^2}{(\text{trace}(I_k - \omega K C({\lambda}_{x}, {\lambda}_{\xi})))^2},
\end{equation}
where $\omega = k/m$ \cite{renaut2017hybrid}, and $C({\lambda}_{x}, {\lambda}_{\xi})$ is defined in \eqref{eq:C}.  Note that other strategies for selecting suitable regularization parameters for each of the projected problems can be seamlessly adopted within this framework.  For example, by generalizing to the two-parameter case other common alternatives typically employed within hybrid methods, see e.g. \cite{hybrid_review,msl2}. 

For the examples shown in Section \ref{Chapter:experiments}, the proposed algorithms are also equipped with a stopping criterion. Following \cite{sdecompose}, we use the flattening of the generalized cross validation (GCV) function defined with respect to the iterations,
\begin{equation} \label{eq:GCVstop}
    \widehat{G}(k) = \frac{k \|K y_k({\lambda}^{(k)}_{x}, {\lambda}^{(k)}_{\xi}) - \|b\|_{R^{-1}} e_1\|_2^2 }{(\text{trace}(I_k -  K C (\lambda^{(k)}_{x},\lambda^{(k)}_{\xi})))^2}.
\end{equation}
If the stopping criterion is not satisfied, the iterations are stopped after a given maximum amount of iterations.

Another algorithmic choice that needs to be mentioned is the initialization of the weights $W_1$, which are constructed using the first reconstruction available when only one basis vector is added in the space, i.e. $W_1 = W(u_1) = W(x_1).$ Lastly, we have empirically observed that the algorithm is very robust to the choice of $\tau$ in \eqref{eq:weights}, as long as it is some orders of magnitude smaller than then average pixel intensity, and bounded away from machine precision. 

\section{Convergence of AF-GMRES and AF-LSQR}\label{Chapter:proof}

\textbf{Lemma 3.1.}  Assume that  no  breakdown  happens  in the augmented flexible Arnoldi and the augmented flexible Golub–Kahan algorithms.  Then, the sequence
$\phi(x_k,\xi_k)$, where $\phi(x, \xi)$ is defined in \eqref{eqn:def_phi} and $x_k,$ $\xi_k$ are the approximate solution computed after $k$ steps of AF-GMRES or AF-LSQR, is decreasing monotonically and it is bounded from below by zero.

\textit{Proof} Since $\phi(x, \xi) \geq 0$, we only need to prove that $\phi(x_k,\xi_k) \leq \phi(x_{k-1},\xi_{k-1})$. To do this recall that $\phi_k(x,\xi)$ is a quadratic tangent majorant of $\phi(x,\xi)$ at the points $x,$ $\xi_{k-1}$ and, in particular, at $x_{k-1},$ $\xi_{k-1}$, so that:
\begin{equation}\label{eqn:proof}
\phi(x,\xi) \leq \phi_k(x,\xi) \quad \forall x,\xi  \in \mathbb{R}^{n}, \forall k
\quad \text{and} \quad \phi_k(x_{k-1},\xi_{k-1}) = \phi(x_{k-1},\xi_{k-1}) \quad \forall k.
\end{equation}
Then, for  $1 \leq k \leq n/2$
\begin{equation}
\phi(x_k,\xi_k) \leq \phi_k(x_k,\xi_k) = \min_{\substack{x \in \mathcal{R}(V_k) \\ \xi \in \mathcal{R}(Z_{k-1})}} \phi_k(x,\xi) \leq \phi_k(x_{k-1},\xi_{k-1}) = \phi(x_{k-1},\xi_{k-1}).
\end{equation}
Here, the first inequality uses the first expression in \eqref{eqn:proof}, the second inequality uses that $x_{k-1} \in \mathcal{R}(V_{k-1}) \subset \mathcal{R}(V_k)$ and $\xi_{k-1} \in \mathcal{R}(Z_{k-2}) \subset \mathcal{R}(Z_{k-1})$, and the last equality holds because of the second expression in \eqref{eqn:proof}. Finally, for  $k \geq n/2$ 
\begin{equation}
\phi(x_k,\xi_k) \leq \phi_k(x_k,\xi_k) = \min_{\substack{x \in \mathbb{R}^{n} \\ \xi \in \mathbb{R}^{n}} } \phi_k(x,\xi) \leq \phi_k(x_{k-1},\xi_{k-1}) = \phi(x_{k-1},\xi_{k-1}).
\end{equation}\\ 
\textbf{Theorem 3.2.} Under the same assumptions of Lemma 3.1. the sequence $\{x_k,\xi_k\}_{k\geq1}$, where $x_k$, $\xi_k$ is the approximate solution computed after $k$ steps of AF-GMRES or AF-LSQR is such that $\lim_{\substack{k \rightarrow \infty}} \|x_k-x_{k-1}\|_2 = 0$ and $\lim_{\substack{k \rightarrow \infty}} \|\xi_k-\xi_{k-1}\|_2 = 0$. Moreover, it converges to a stationary point of $\phi(x,\xi)$.

\textit{Proof} Thanks to Lemma 3.1, $\phi(x,\xi) \geq 0$ has a stationary point.  An intuitive way of proving strong convexity of the functionals $\phi_k(x,\xi)$ for all $k$ is to re-write them in the following form:
\begin{equation}\label{eq:proof}
    \bar{\phi}_k( \begin{bmatrix}
        x \\ \xi
    \end{bmatrix}) = \left|\left| R^{-1/2}\left(A  \begin{bmatrix}
        I & I \end{bmatrix}\right)
        \begin{bmatrix}
        x \\ \xi
    \end{bmatrix} - b \right|\right|_2^2 + \left|\left| 
    \begin{bmatrix}
    {\lambda}_{x} Q^{1/2} & 0 \\
0  & {\lambda}_{\xi} {W_{k-1}} 
\end{bmatrix}  \begin{bmatrix}
        x \\ \xi
\end{bmatrix}\right|\right|_2^2.
\end{equation}
Now it can be clearly seen that the functional in \eqref{eq:proof} is strongly convex, as the first term is convex and the second term is strongly convex (its Hessian has strictly positive eigenvalues). Note that \eqref{eq:proof} is never constructed in practice. However, for $k \geq n/2$, $[x_k,\, \xi_{k}]^{\top}$ is the minimizer of \eqref{eq:proof}, and by strong convexity of \eqref{eq:proof};
\begin{equation}
    \frac{\rho}{2}  \left\|\begin{bmatrix}
        x_k \\ \xi_k
    \end{bmatrix} - \begin{bmatrix}
        x_{k-1} \\ \xi_{k-1}
    \end{bmatrix}\right\| \leq \bar{\phi}_{k}(\begin{bmatrix}
        x_{k-1} \\ \xi_{k-1}
    \end{bmatrix})- \bar{\phi}_{k}(\begin{bmatrix}
        x_{k} \\ \xi_{k}
    \end{bmatrix}) \leq \bar{\phi}(\begin{bmatrix}
        x_{k-1} \\ \xi_{k-1}
    \end{bmatrix})- \bar{\phi}(\begin{bmatrix}
        x_k \\ \xi_k
    \end{bmatrix}) \underset{k\rightarrow\infty}{\rightarrow} 0
\end{equation}
where we have used that $\phi_k$ is an upper bound of $\phi$ and that $\phi(x_{k-1},{\xi}_{k-1})- \phi(x_k, \xi_k) \underset{k\rightarrow\infty}{\rightarrow} 0$ by Lemma 3.1. Therefore, $\lim_{\substack{k \rightarrow \infty}} \|x_k-x_{k-1}\|_2 = 0$ and $\lim_{\substack{k \rightarrow \infty}} \|\xi_k-\xi_{k-1}\|_2 = 0$.

\section{Results}\label{Chapter:experiments} 
In this section, three experiments are presented to investigate the performance of the proposed methods. In Subsection \ref{subsec:deblur}, we show an image deblurring problem with a square system matrix to showcase the performance of AF-GMRES in comparison to hybrid GMRES methods and hybrid flexible GMRES (denoted by FGMRES) methods. In the subsequent two subsections, we evaluate the efficacy of AF-LSQR in comparison to hybrid flexible LSQR (FLSQR) that only incorporates a sparsity-inducing $\ell_1$ regularization term, the generalized hybrid method (genHyBR) that exclusively uses a $Q$-weighted regularization term, and the solution decomposition hybrid method (sdHyBR) that employs a flexible generalized Golub-Kahan iterative approach to solve \eqref{eq:std_form}. In Subsection \ref{subsec:recon}, we consider a synthetic atmospheric transport problem, where the solution corresponds to modeled gas emissions in North America. In Subsection \ref{subsec:fmt}, we present a phantom fluorescence molecular tomography experiment to demonstrate the effectiveness of AF-LSQR with both simulated and real-world data. 

In the first two examples, the discrepancy principle (DP) was employed to select all the relevant regularization parameters in each minimization. In the last experiment,
WGCV was used instead because the noise level for this experiment was unknown. In AF-GMRES, AF-LSQR and sdHyBR, choosing the regularization parameters requires solving two-dimensional nonlinear constraint optimization problems, where we use a Quasi-Newton method as implemented in MATLAB's fminunc function with an initial guess of $\lambda_x = -0.5$ and $\lambda_{\xi} = -0.5$. For the stopping criteria, the iterative process is terminated if either of the following two conditions is satisfied: (i) a maximum number of iterations is reached, or (ii) the GCV stopping function defined in \eqref{eq:GCVstop} flattens out. The following experiments were run on a laptop computer with Intel i5 CPU 2GHz and 16G memory.

\subsection{Image deblurring example}\label{subsec:deblur}
\begin{figure}[ht!]
        \centering
    \begin{tabular}{ccc}
       \includegraphics[width = 0.3\textwidth]{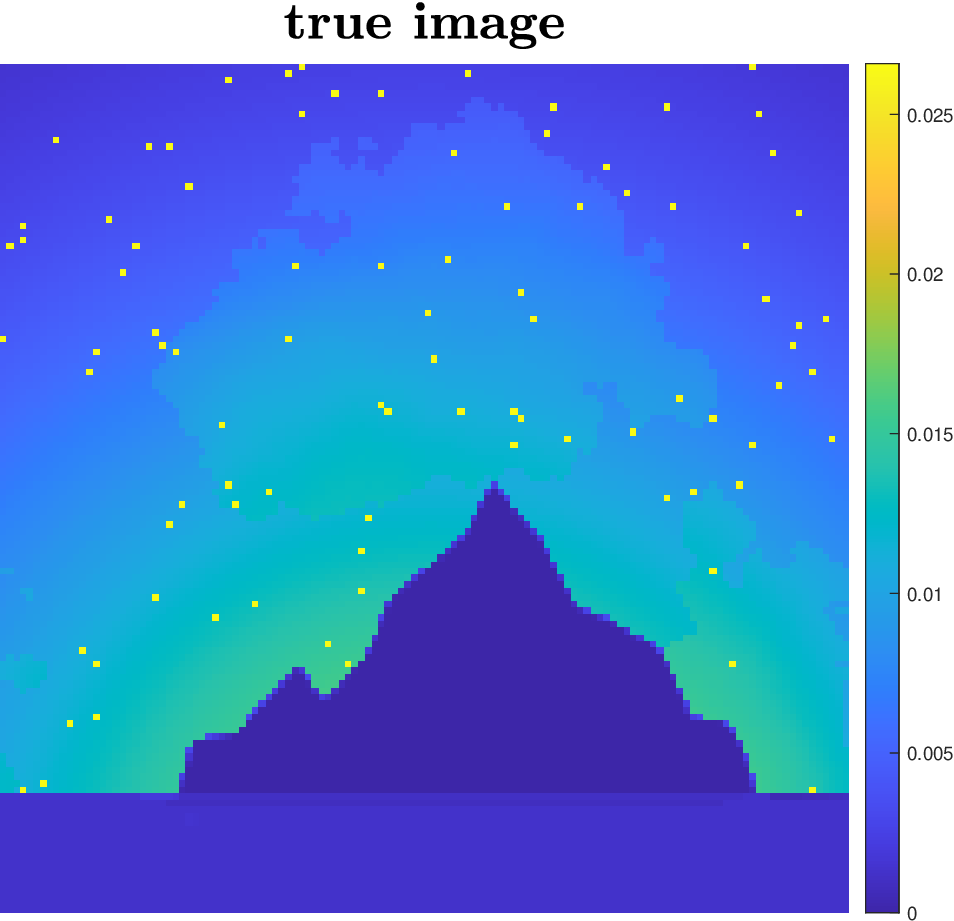}  & \includegraphics[width = 0.3\textwidth]{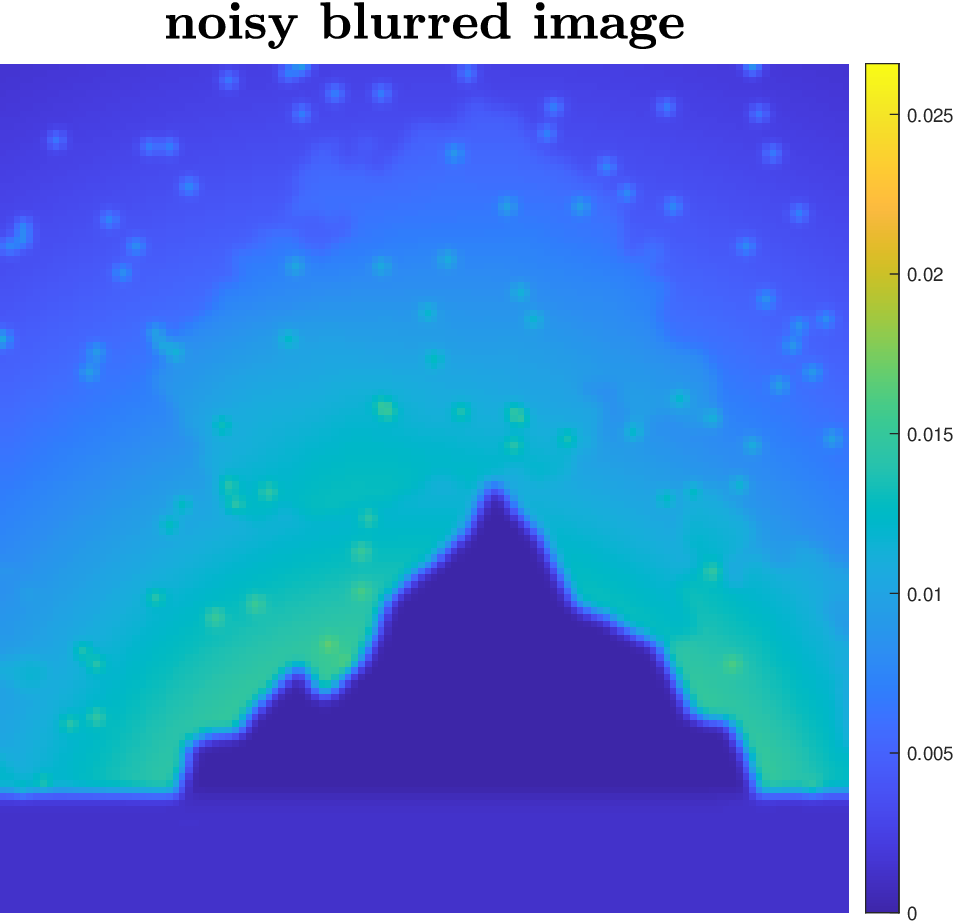} 
         & \includegraphics[width = 0.3\textwidth]{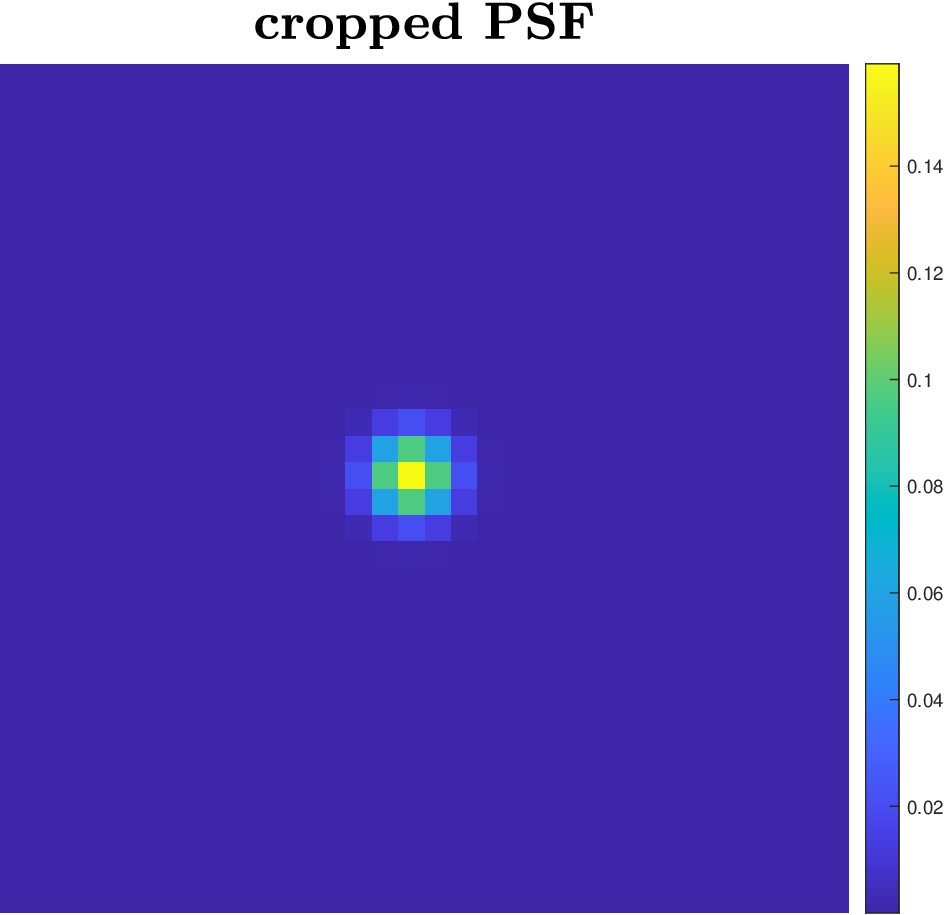}\\
    \end{tabular}
    \caption{Image deblurring example}
    \label{fig1_ex1}
\end{figure}

The  first  experiment  corresponds to an 
image deblurring problem where the exact test  image  of  size  $128 \times 128$  pixels has been chosen to contain both smooth and sparse features. The system matrix representing Gaussian blurring with variance $\sigma =1$, has been computed using IRtools \cite{gazzola2019ir}, and the measurements are corrupted by Gaussian white noise of level $\frac{\|e\|_2}{\|Au_{\text{true}}\|_2} = 10^{-5}$. The true image, blurred and noisy image, and point spread function (PSF) for this example are shown in  Figure \ref{fig1_ex1}.

\begin{figure}[ht!]
    \centering
    \begin{tabular}{ccc}
    \includegraphics[width = 0.4\textwidth]{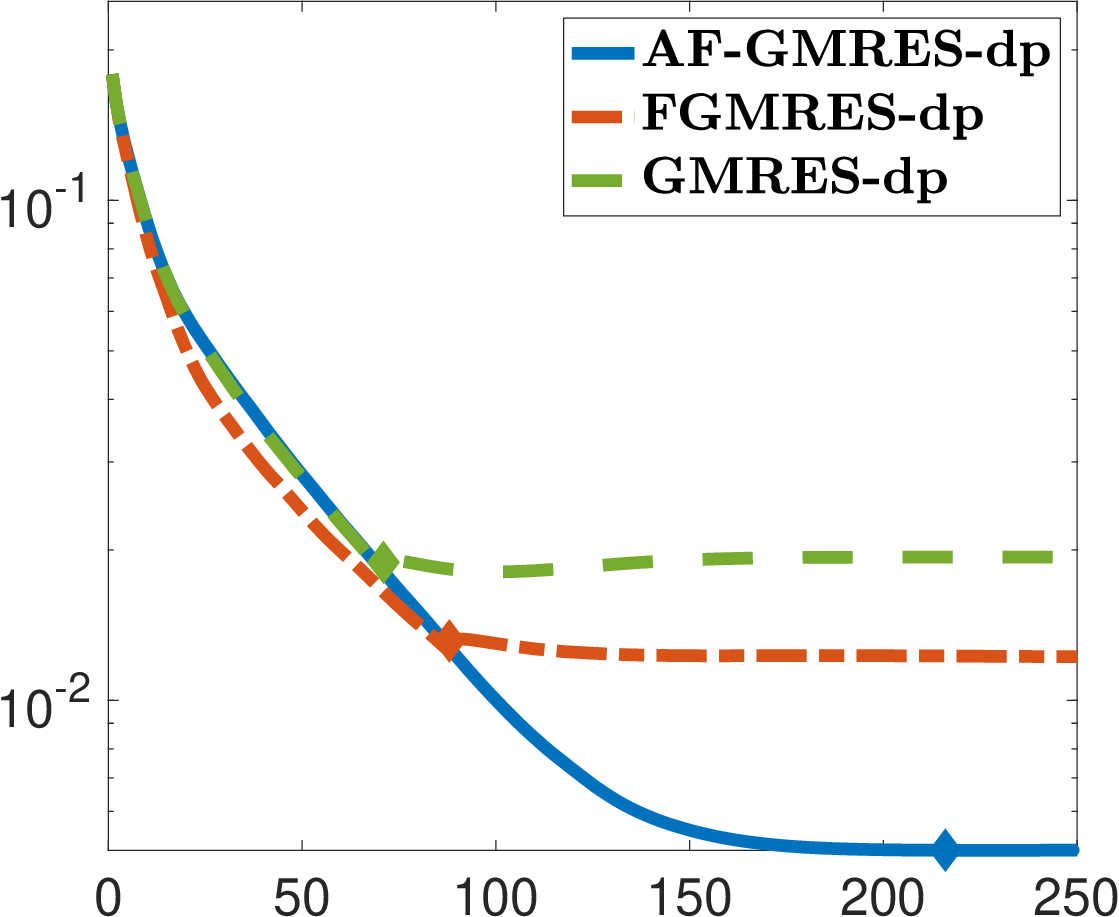} \hspace{0.5cm}
    &     \includegraphics[width = 0.4\textwidth]{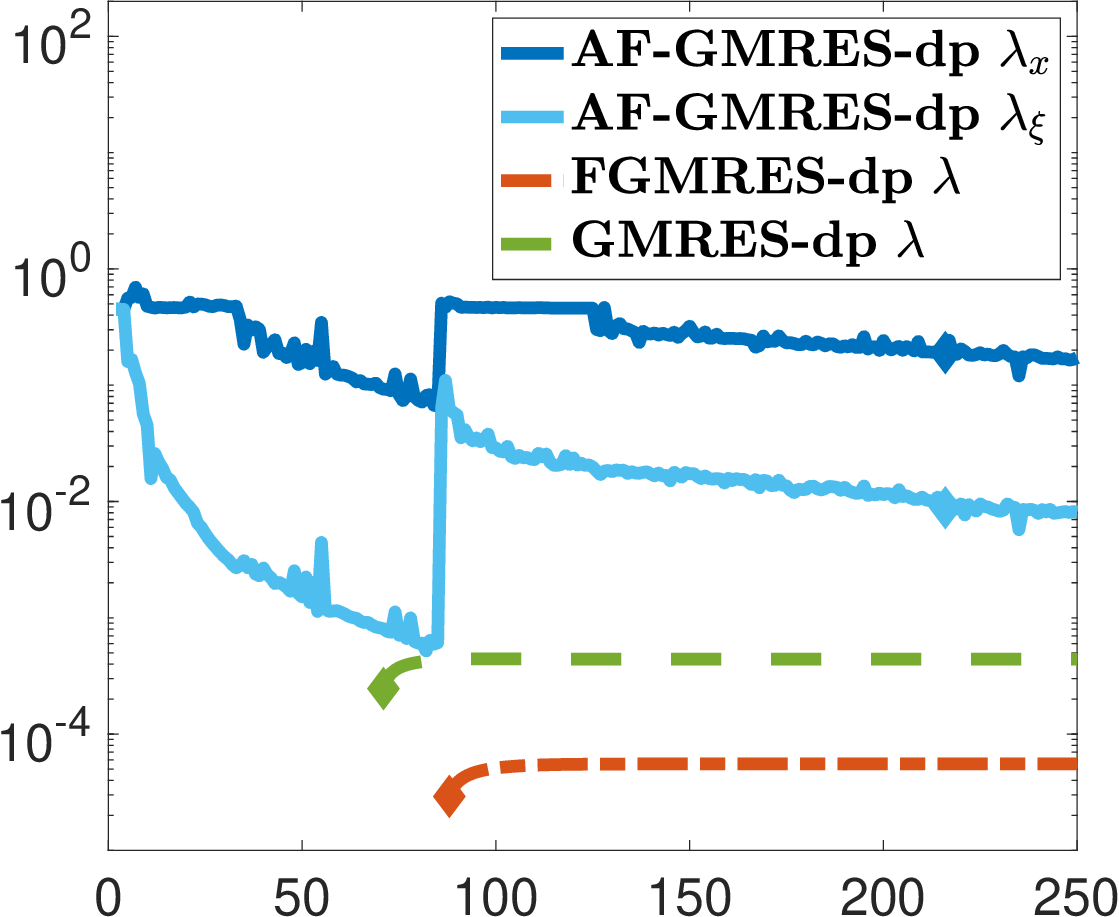} \\
         (a) Relative error norms & (b) Regularization parameters
    \end{tabular}
    \caption{(a) Relative error norms and (b) regularization parameters per iteration for the problem in Figure \ref{fig1_ex1} for AF-GMRES and other GMRES-based solvers. The regularization parameters are computed at each iteration using the DP in equation \eqref{eq:dp}, with $\tau_{dp}=1.1$. The flattening of the GCV function in \eqref{eq:GCVstop} is used to stop the iterations, with a tolerance of $0.02$, and the stopping iteration is signaled in the plots with a marker. }
    \label{fig2_ex1}
\end{figure}

Since $A$ is square, the performance of AF-GMRES can be compared to other GMRES-based standard solvers. It can be observed in Figure \ref{fig2_ex1}(a) that AF-GMRES outperforms the compared methods in terms of the relative error norm, where the diamond-shaped markers denote the stopping iterations. The regularization parameters $\lambda^{(k)}_x$ and $\lambda^{(k)}_\xi$ are chosen at each iteration using the DP with $\tau_{dp}=1.1$. It is also interesting to note in Figure \ref{fig2_ex1}(b) that this criterion leads to the stabilization of the regularization parameters. For this example, the flattening of the GCV function in \eqref{eq:GCVstop} is used to stop the iterations with a tolerance of $0.02$.

\begin{figure}[ht!]
    \begin{tabular}{ccc}
    \centering
       \includegraphics[width = 0.3\textwidth]{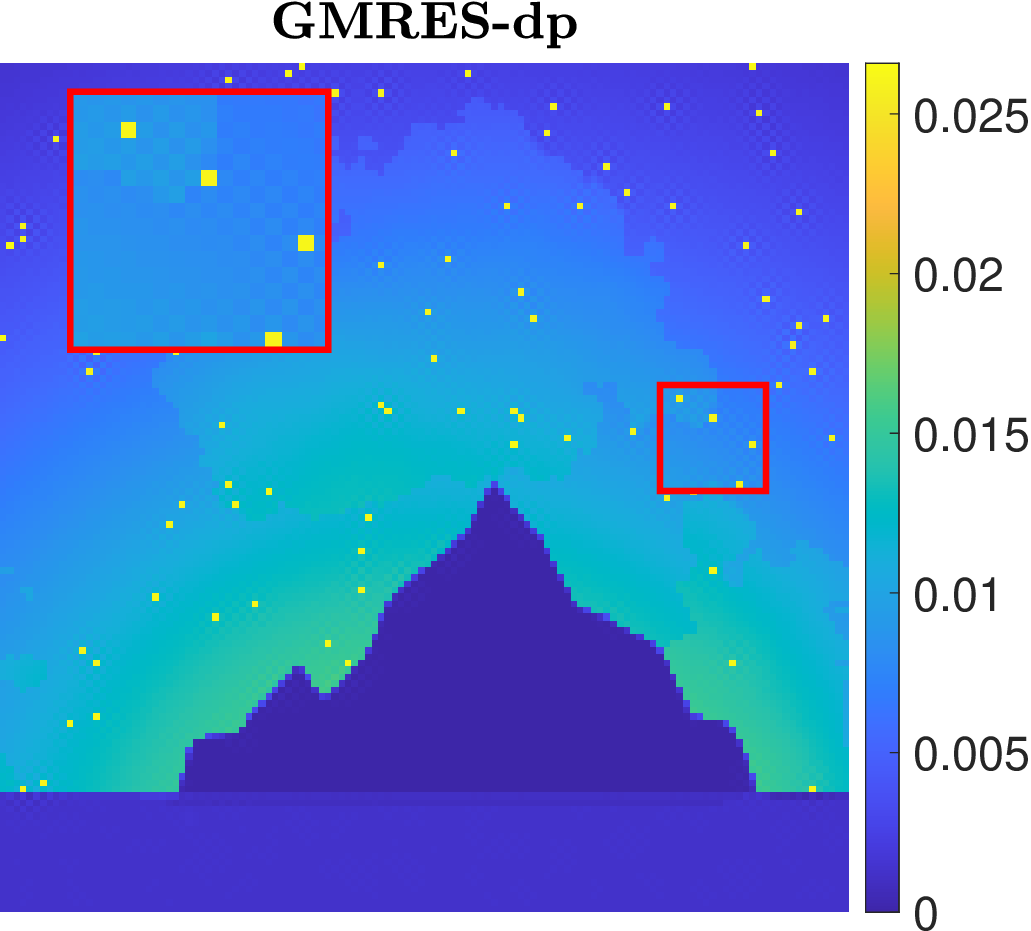}  & \includegraphics[width = 0.3\textwidth]{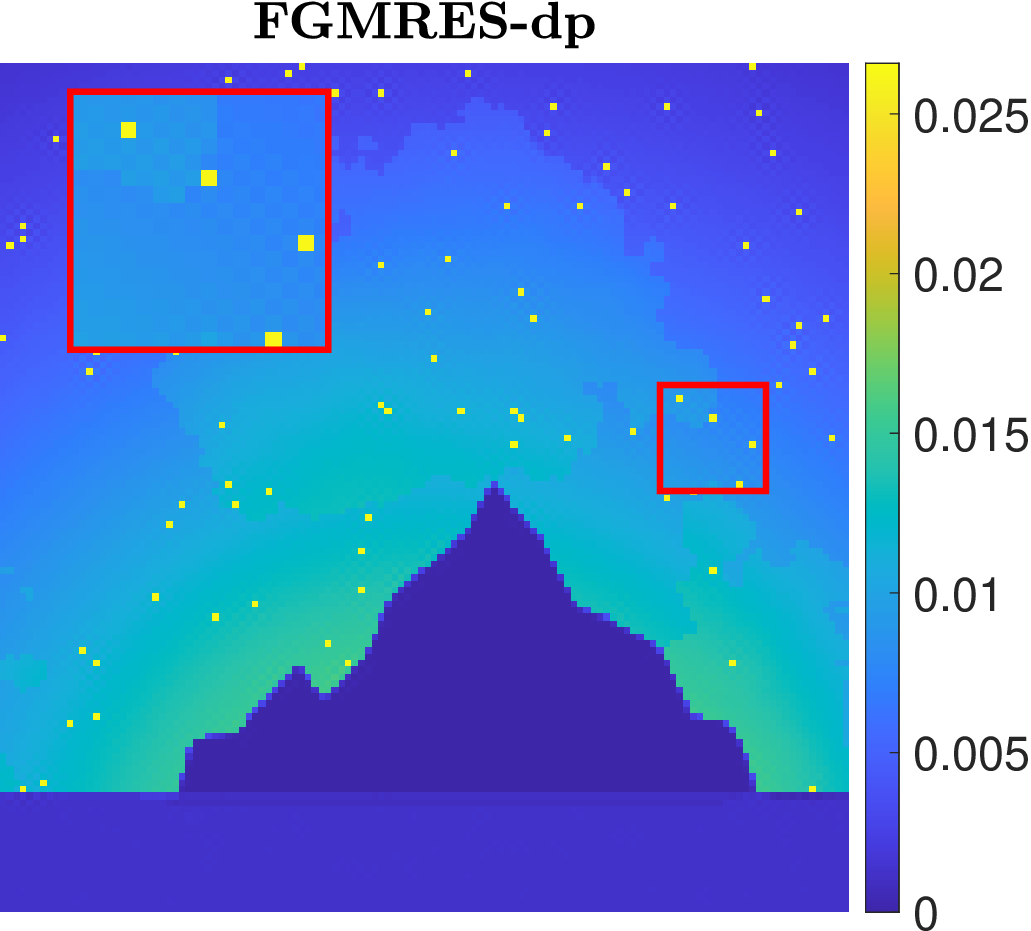} 
         & \includegraphics
         [width = 0.3\textwidth]{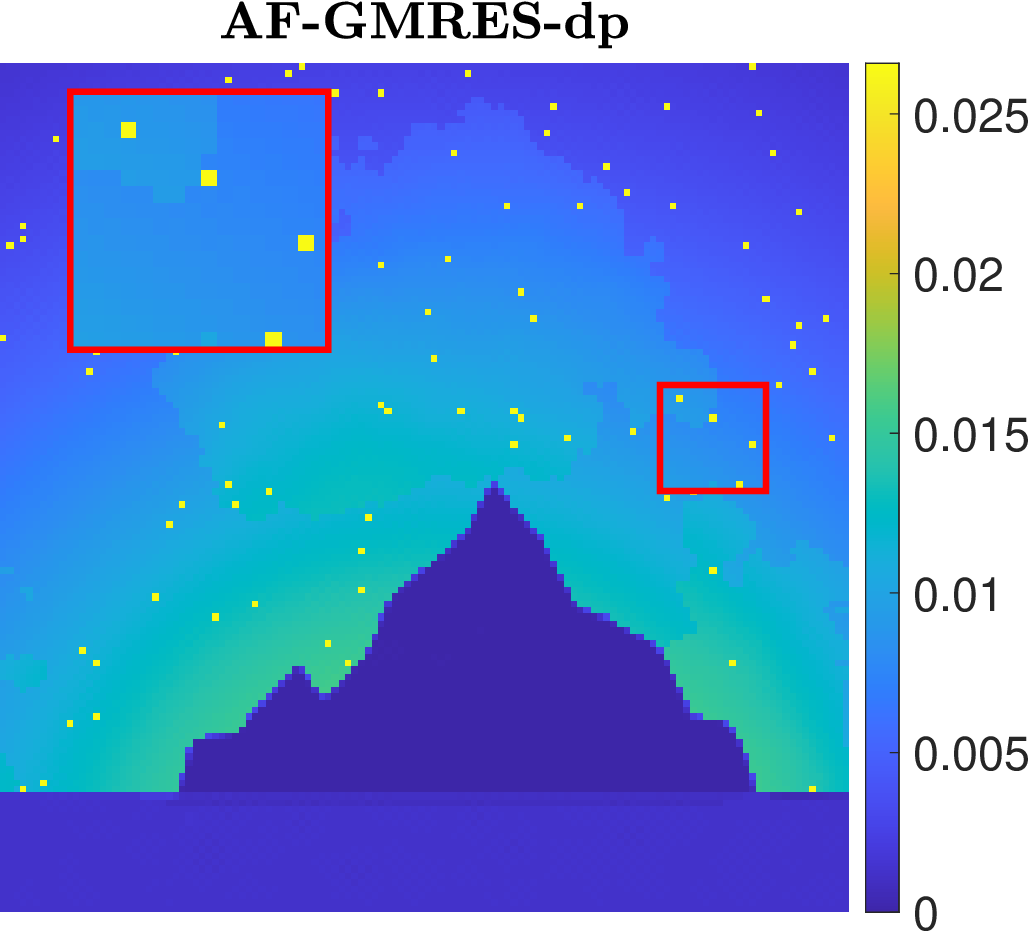} \\
        \end{tabular}
    \caption{Reconstructions at the stopping iteration using the different GMRES-based methods used for comparison. On the top left corner, a zoomed version of the highlighted region in the image is shown to illustrate the details of the reconstruction. All results correspond to using the DP to select the regularization parameters.}
    \label{fig3_ex1}
\end{figure}

It is worth mentioning that although AF-GMRES converges slower than other methods, its reconstructions display a much better quality. This is not only supported by the relative error norm, but can also be directly observed in Figure \ref{fig3_ex1}, where close-ups on the detailed structure are shown in the top left corners.

\subsection{Synthetic atmospheric transport problem.}\label{subsec:recon} 
\begin{figure}[ht!]
\centering
   \hspace{-0.5in}
    \begin{tabular}{cc}
       \includegraphics[height = 6 cm]{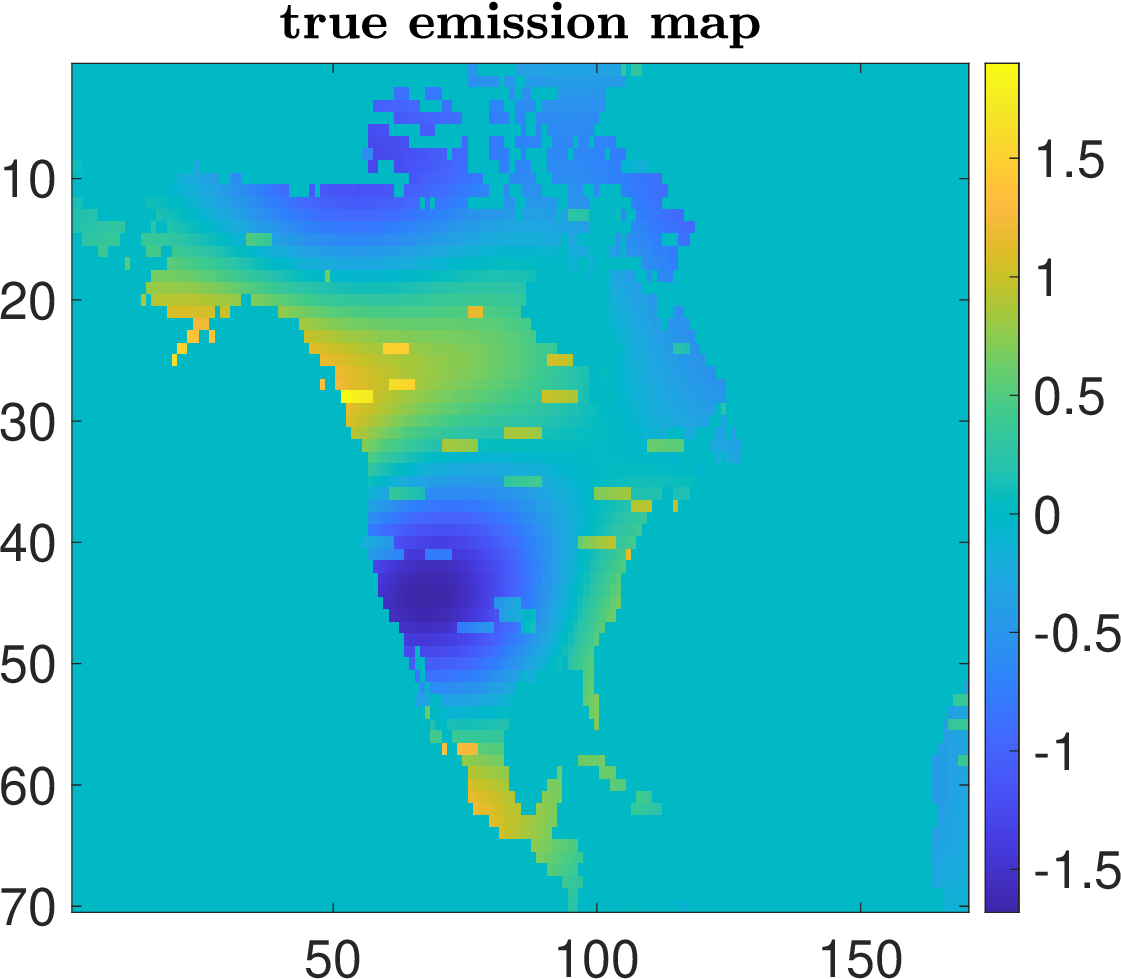} \hspace{0.5cm} & \includegraphics[height = 6 cm]{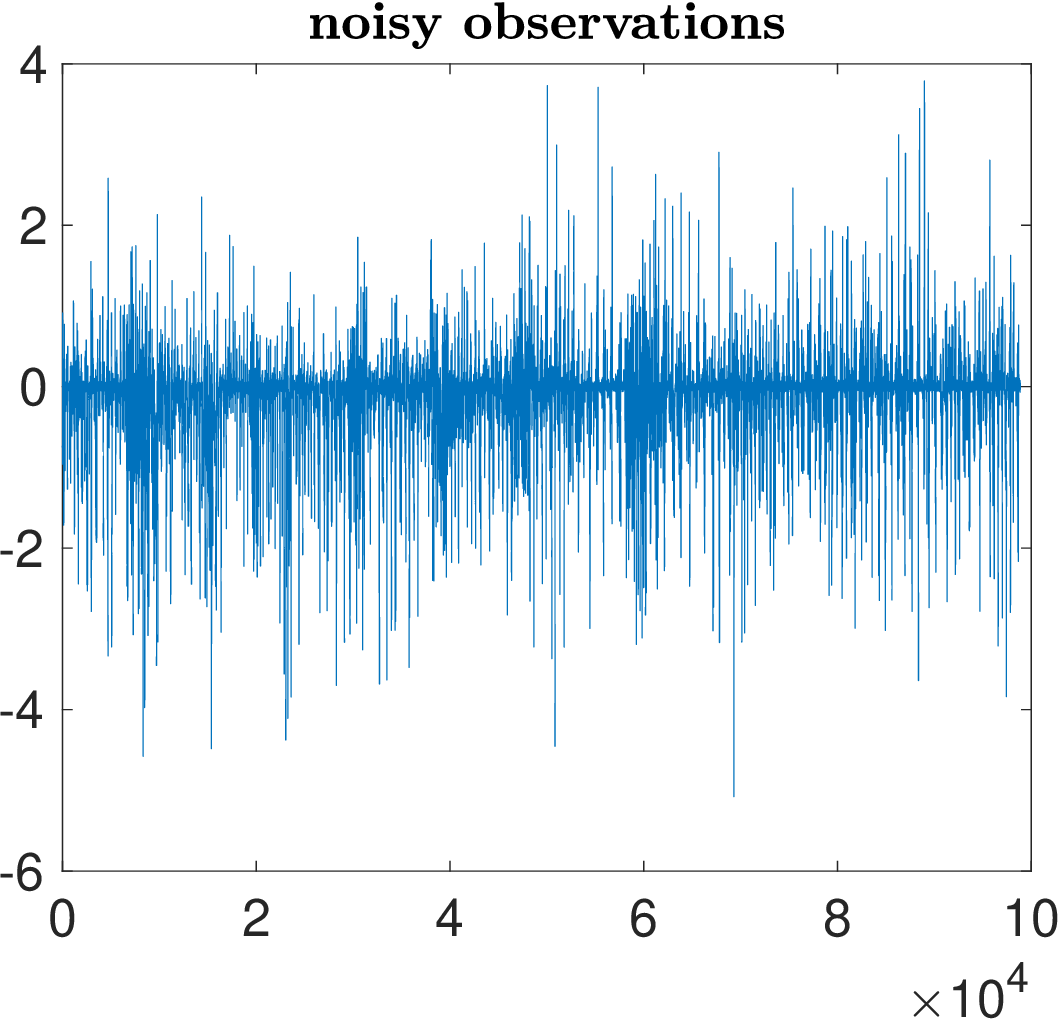}
    \end{tabular}
    \caption{Atmospheric transport example}
    \label{fig1_ex2}
\end{figure}

In this experiment, we consider a synthetic atmospheric transport problem reconstructing a carbon dioxide (CO$_2$) flux emission map that includes sparse high CO$_2$ fluxes and a smooth CO$_2$ flux background.  Figure \ref{fig1_ex2} presents the true emissions, derived from the sum of a smooth image generated using a Matérn kernel \cite[Equation (4.14)]{10.7551/mitpress/3206.001.0001} with parameters $\nu = 2.5$ and $\ell = 0.05$, and an image with sparse speckles whose intensity matches the maximum value of the smooth image. The forward atmospheric model $A \in \mathbb{R}^{98880 \times 3222}$ is taken from the NOAA's Carbon Tracker Lagrange project \cite{liu2021data,miller2020geostatistical} and was created using the stochastic time-Inverted Lagrangian transport system \cite{nehrkorn2010coupled,lin2003near}. The observations $b \in \mathbb{R}^{98880}$ shown in Figure \ref{fig1_ex2} are measurements at the locations and times measured by the Orbiting Carbon Observatory 2 (OCO-2) satellite between July and mid-August 2015, with added Gaussian white noise with noise level of $\frac{\|e\|_2}{\|Au_{\text{true}}\|_2} = 0.04$. The spatial resolution modeled in this example is $1^{\circ} \times 1^{\circ}$. Although this resolution is not fine enough to detect realistic super-emitters in practice, these examples are provided as a proof of concept for the studied methods, with the note that these can be extended to larger datasets.  To avoid committing the inverse crime, the covariance matrix $Q$ used in the functional minimized by genHyBR and AF-LSQR is taken to be a Matérn kernel with $\nu = 1$ and $\ell = 0.5$. 
The flattening of the GCV function is used as a stopping criterion with a tolerance of  $5\times10^{-6}, 10^{-7}, 10^{-7}, 10^{-6}$ for AF-LSQR, sdHyBR, genHyBR and FLSQR respectively. 

\begin{figure}[ht!]
    \centering
    \includegraphics[width = 0.4\textwidth]{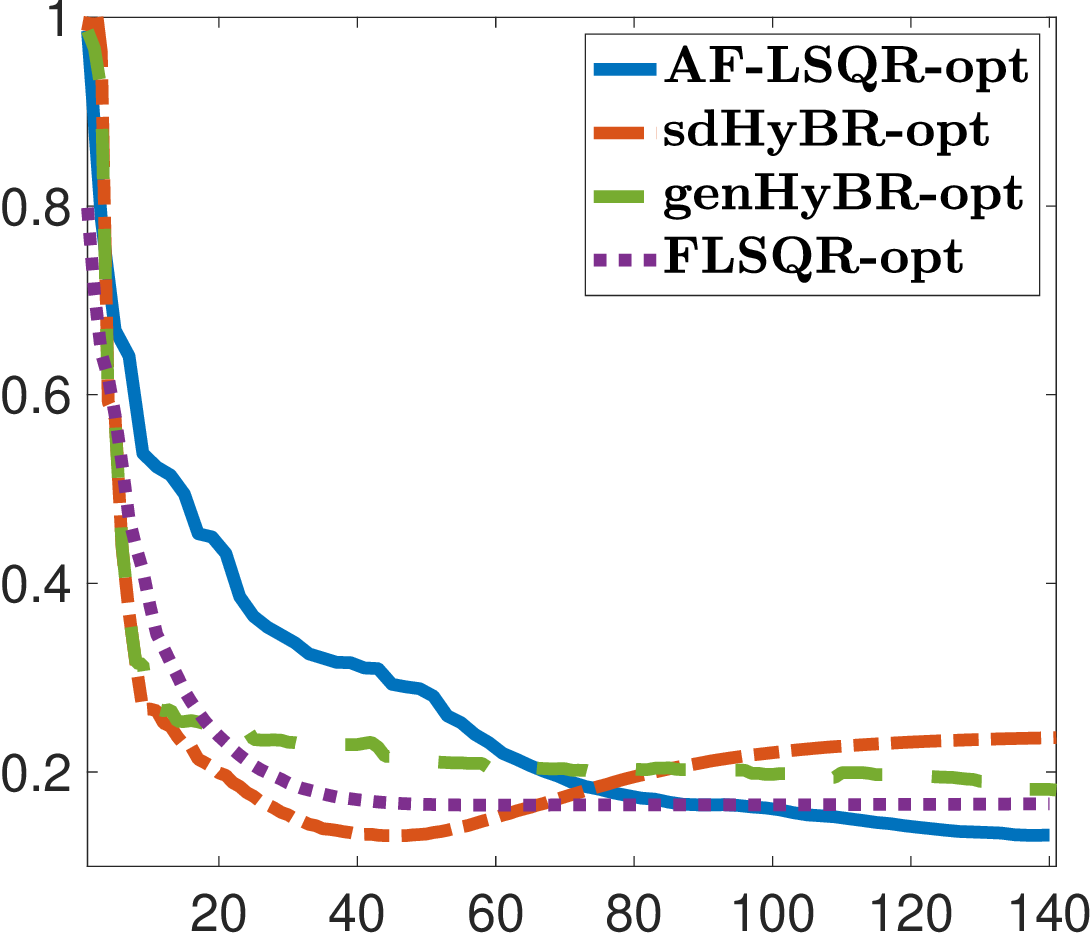}
    \hspace{0.5cm} \includegraphics[width = 0.4\textwidth]{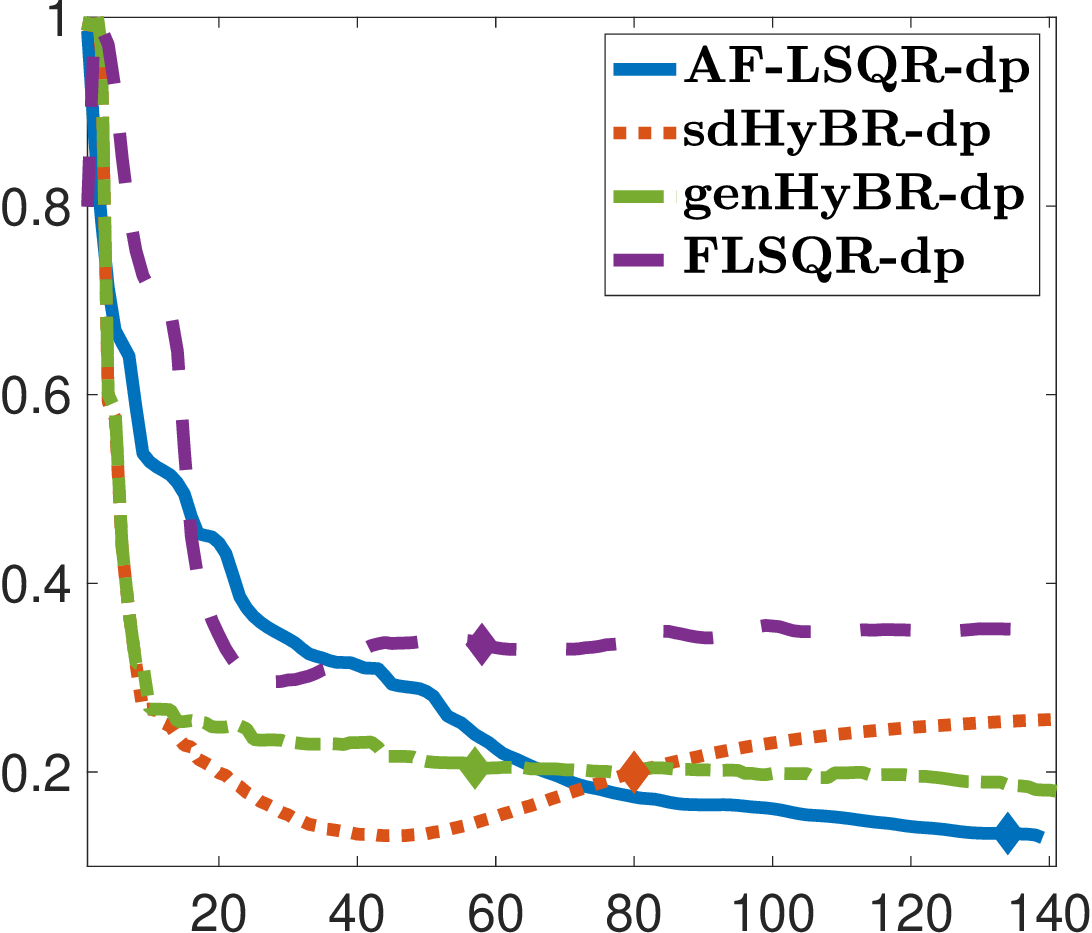}
    \caption{Atmospheric transport example: Relative reconstruction error norms per iteration of sdHyBR, genHyBR, FLSQR and AF-LSQR. Results in the left plot correspond to selecting the optimal regularization parameters at each iteration, and results in the right plot correspond to the DP-selected regularization parameters.}
    \label{fig2_ex2}
\end{figure}

\begin{figure}[ht!]
\center
    \begin{tabular}{cc}
    \centering
       \includegraphics[width = 0.4\textwidth]{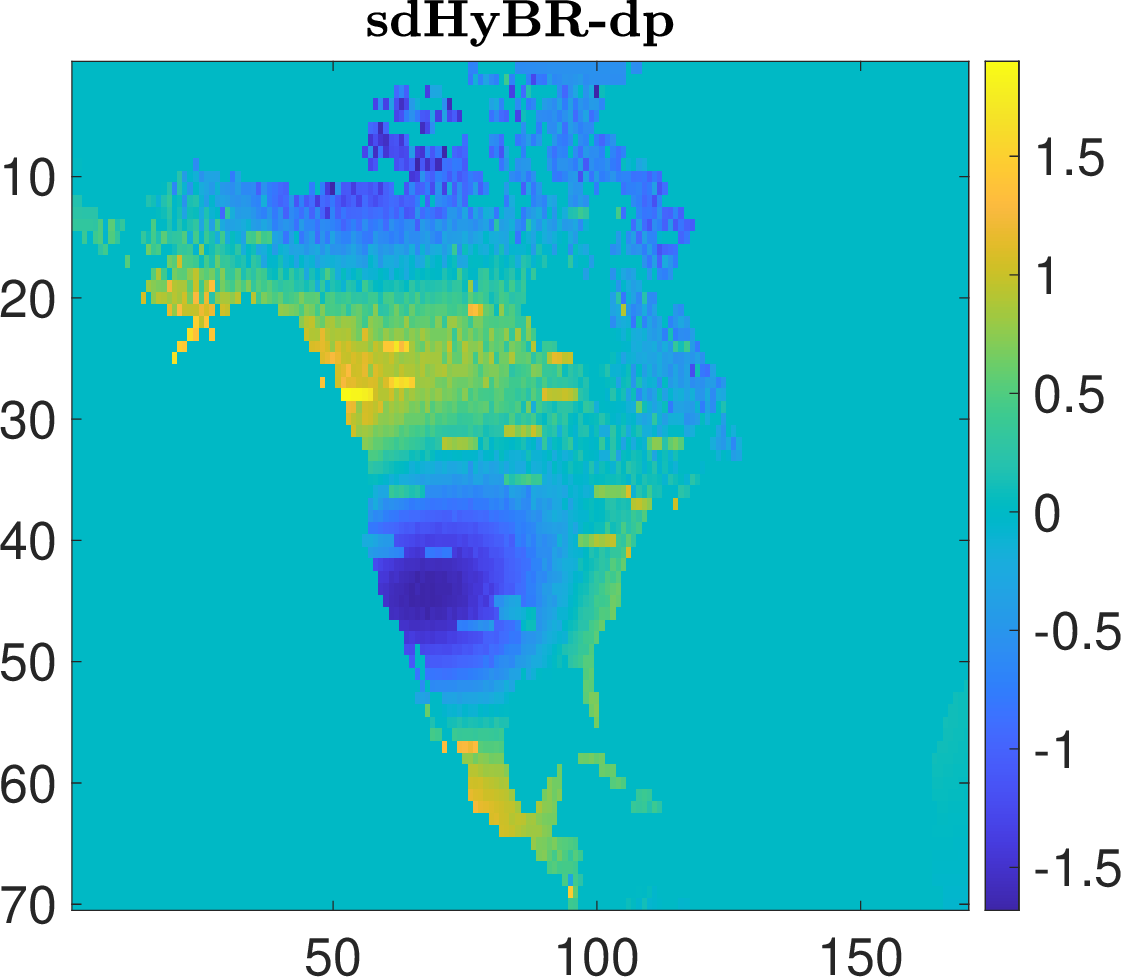}  \hspace{0.5cm} & \includegraphics[width = 0.4\textwidth]{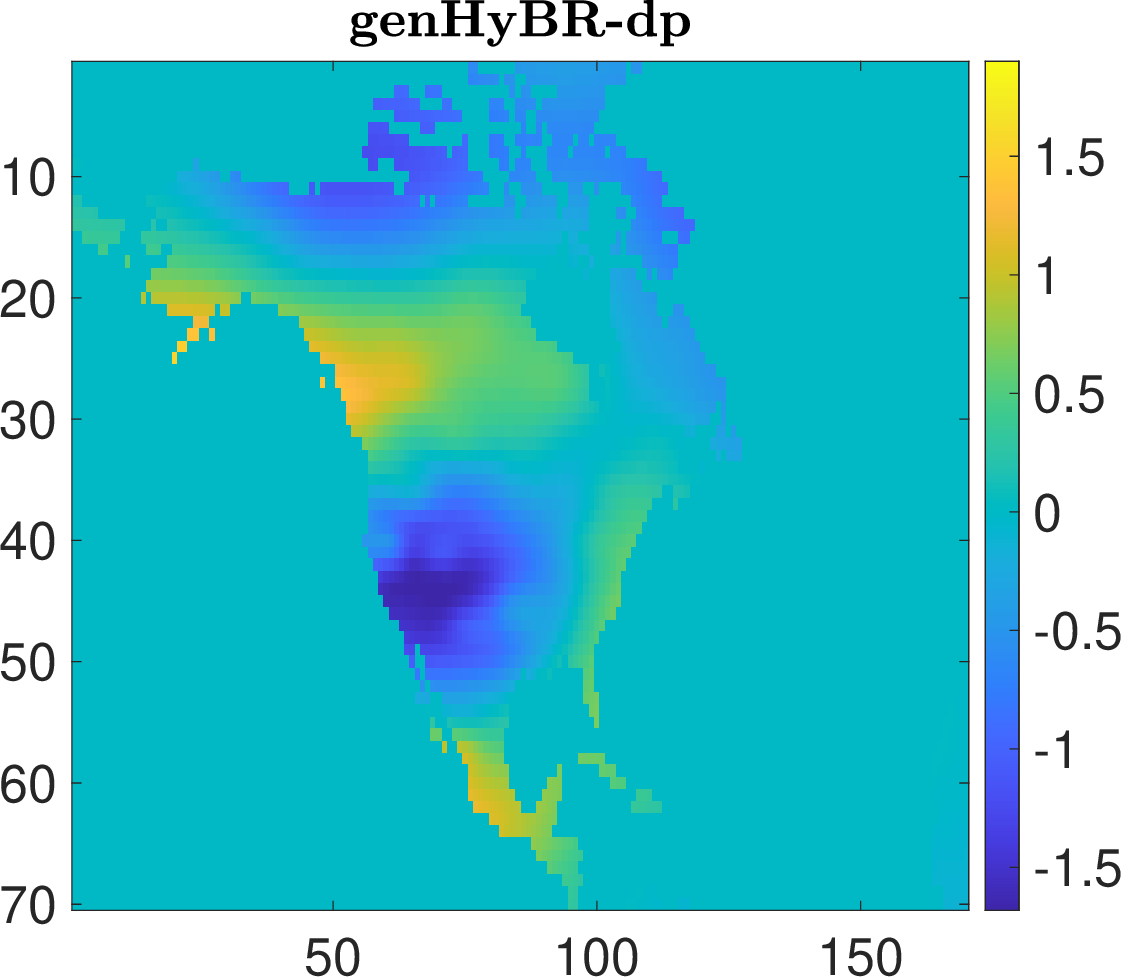} \\
         \includegraphics[width = 0.4\textwidth]{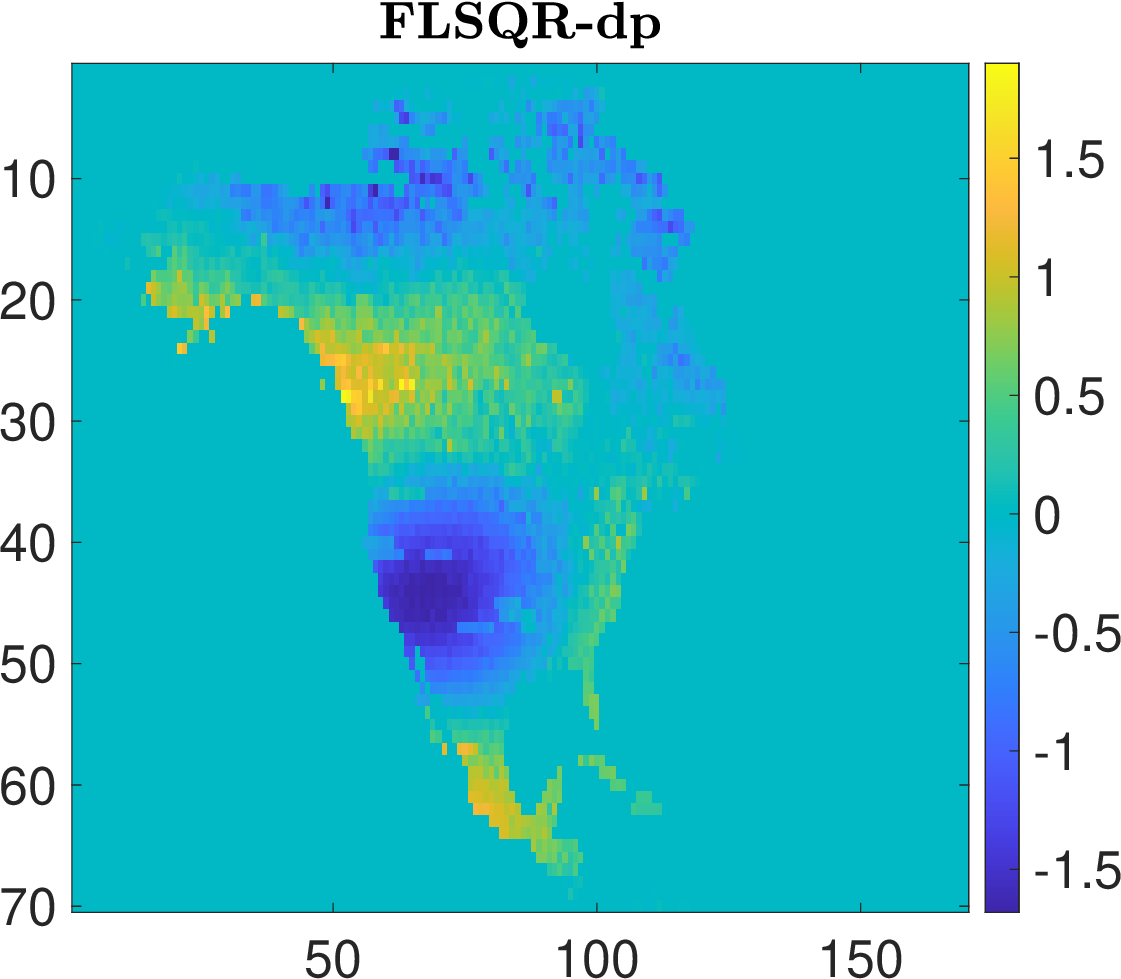}  \hspace{0.5cm} & \includegraphics[width = 0.4\textwidth]{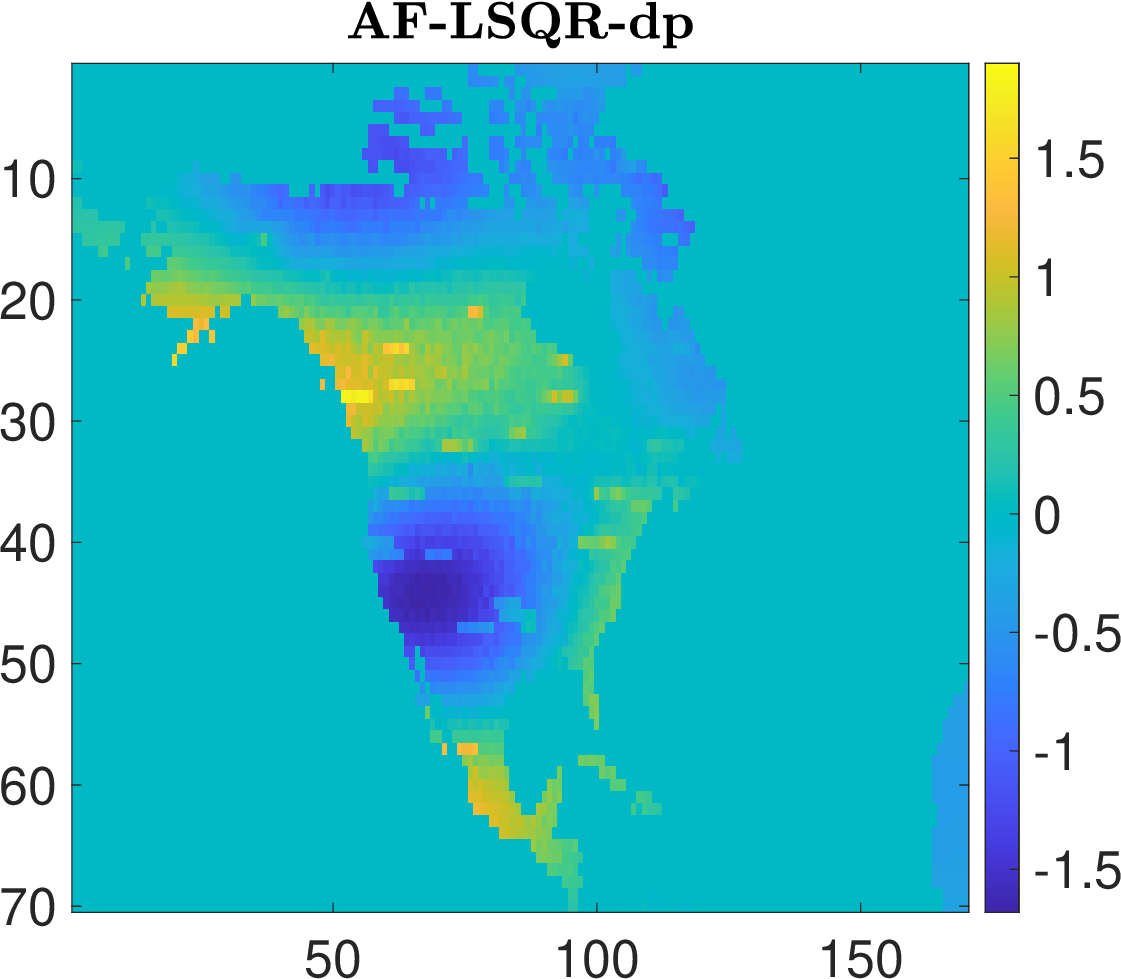}
    \end{tabular}
    \caption{Atmospheric transport example: In the top row, we provide the reconstruction of the fluxes for sdHyBR and genHyBR. In the bottom row, we provide the reconstruction of the fluxes for FLSQR and AF-LSQR. All results correspond to using DP-selected regularization parameters and are obtained at the stopping iterations.}
    \label{fig1_ex2_rec}
\end{figure}

In the left plot of Figure \eqref{fig2_ex2}, relative reconstruction error norms per iteration are provided using the optimal regularization parameter. In the right plot of Figure \eqref{fig2_ex2}, one can observe that similar results are obtained when using DP-selected regularization parameters. Here, the diamond-shaped markers denote the stopping iterations. We observe that AF-LSQR yields smaller relative reconstruction error than other methods. This result is further illustrated by the image reconstructions displayed in Figure \ref{fig1_ex2_rec}.  As expected,  genHyBR effectively captures smooth region but struggles to reconstruct sparse anomalies, since it only includes a weighted $\ell_2$ regularization term. Contrarily, FLSQR successfully captures the sparse features but corrupts the smooth background. Compared with sdHyBR,  in this example, AF-LSQR achieves a better recovery of the smoothness in the background.

Given the general semi-convergence behavior of ill-posed problems, attaining a good reconstruction can be a very difficult when the regularization parameter is not chosen appropriately. In contrast, the proposed AF-LSQR method shows better convergence and robustness, prompting us to further explore its application on a more realistic problem in the following subsection. 

\subsection{Fluorescence molecular tomography (FMT) reconstruction using synthetic and experimental data} \label{subsec:fmt}
Lastly, the proposed AF-LSQR method is evaluated in a medical imaging problem involving synthetic and experimental data. The imaging modality we choose here is called fluorescence molecular tomography (FMT) which is a non-invasive approah to visualize the 3D distribution of fluorescence in biological tissues. FMT has emerged as an indispensable tool for longitudinal and large-scale monitoring of biological processes and disease progression in living animals with a major application in drug development \cite{hilderbrand2010near,ntziachristos2006fluorescence}. The underlying principle of FMT is based on diffuse light propagation modeling within a scattering volume, which leads to high-similarity adjacent projection images and incomplete measured optical signals, making FMT reconstruction a large-scale and highly ill-posed inverse problem \cite{arridge2009optical}. In FMT reconstruction, the meaningful fluorescence signal is typically localized and concentrated in a restricted region such as a tumor, characterized by its sparsity, whereas this sparse and high-intensity signal is normally surrounded by a smoothly changing background associated to healthy tissues. We apply our method in two FMT reconstruction tasks using synthetic and experimental data respectively. Both cases utilize an identical problem setting with a slab phantom containing a sparse fluorescence target and a large low-concentration target. 

\subsubsection{FMT reconstruction using synthetic data} 
\begin{figure}[ht!]
   \includegraphics[width = 0.9\textwidth]{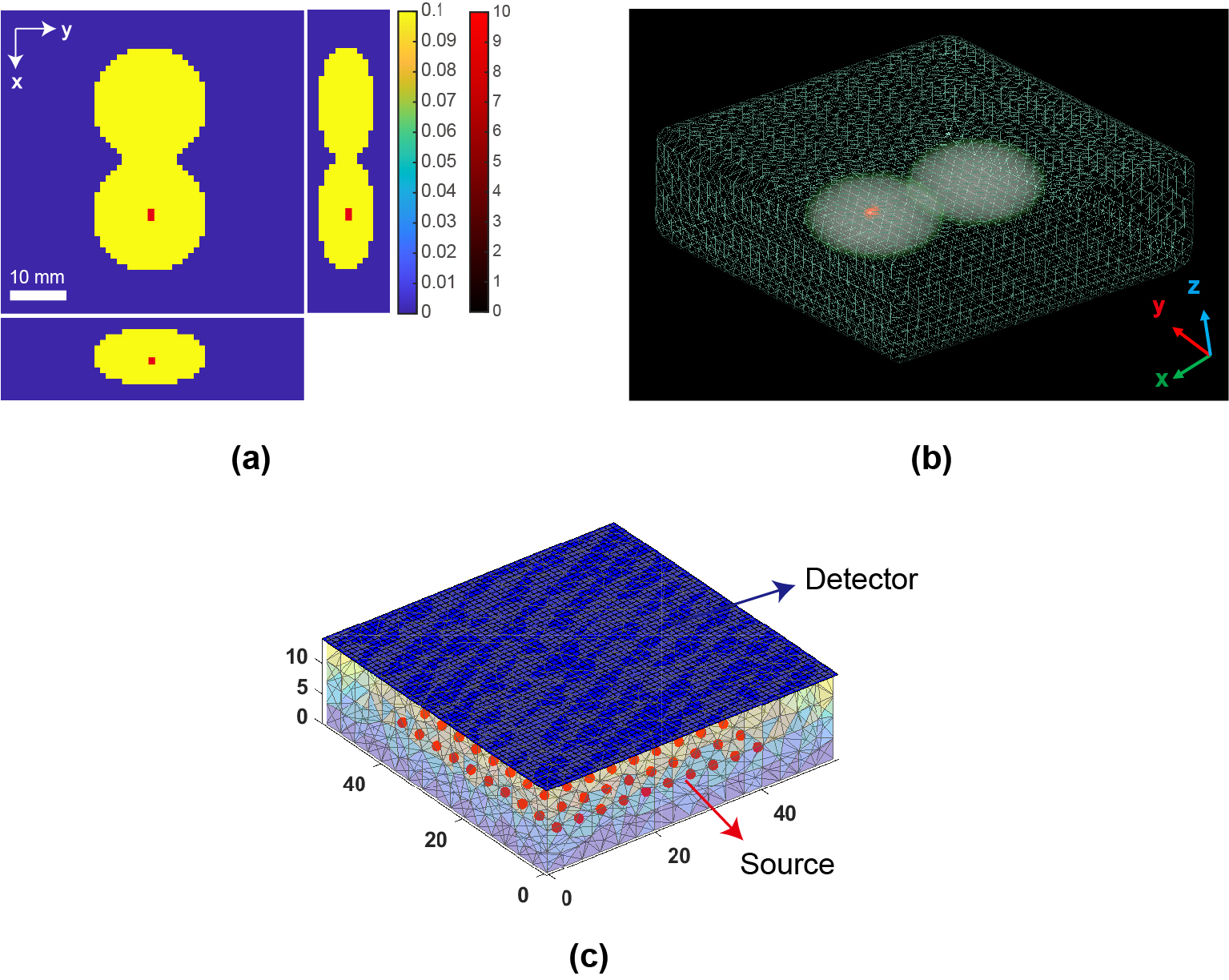}
   \centering
    \caption{Simulation study of FMT: (a) True fluorescence distribution overlayed from smooth and sparse component for the FMT simulation study. Three cross-sections are visualized at the cross position of (38, 28, 8) mm. (b) True fluorescence distribution in a 3D view. The smooth feature and the sparse feature are indicated in green and red respectively. (c) Problem setting of FMT experiment using a simulated slab phantom. A $55\times55$ detector array on the top is simulated for detecting the photon intensity on the surface (blue patches), and $10\times10$ laser points on the bottom surface are applied for illumination (red dots).}
    \label{fig1_ex3}
\end{figure}

In this case study, we consider a simulated phantom model, where the goal is to visualize morphological information of spatially sparse tumors while preserving fluorescence information from the smooth background tissues to aid tumor localization \cite{mieog2022fundamentals}. As depicted in Figure \ref{fig1_ex3}(a-b), the simulated phantom is designed as a slab ($54\times 54 \times 14 $ $mm^3$) and contains two fluorescence inclusions: the first is a large low-concentration region, resembling a peanut, representing the smooth background component, while the second is a small-volume high-concentration component ($2\times 1 \times 1 $ $mm^3$), symbolizing the sparse component.

The simulated transmission-mode detection corresponds to illuminating the bottom of the phantom using a grid of $10\times10$ laser points and collecting fluorescence observations from the opposite side ordered as an array of $55\times55$ covering the whole imaged area as illustrated in Figure \ref{fig1_ex3}(c). The forward matrix $A \in \mathbb{R}^{79577 \times 45375}$ is calculated using the STIFT platform \cite{ren2019smart} and white Gaussian noise with 5$\%$ noise level ($\frac{\|e\|_2}{\|Au_{\text{true}}\|_2} = 0.05$) has been added to the observed fluorescence data $b \in \mathbb{R}^{79577}$. Before computing any reconstruction, the measurements $b$ are pre-processed: the observations that are significantly smaller than $2\%$ of the median of $b$ are considered as noise and therefore removed from $b$. For genHyBR and AF-LSQR, the used covariance matrix $Q$ is generated with a Matérn kernel with parameters $\nu = 10^{-6}$ and $\ell = 0.05$.  The tolerance for the GCV function is chosen to be $2 \times 10^{-6},4 \times 10^{-6}, 10^{-7}, 4 \times 10^{-6}$ for AF-LSQR, sdHyBR, genHyBR and FLSQR respectively. 

\begin{figure}[ht!]
   \includegraphics[width = 0.8\textwidth]{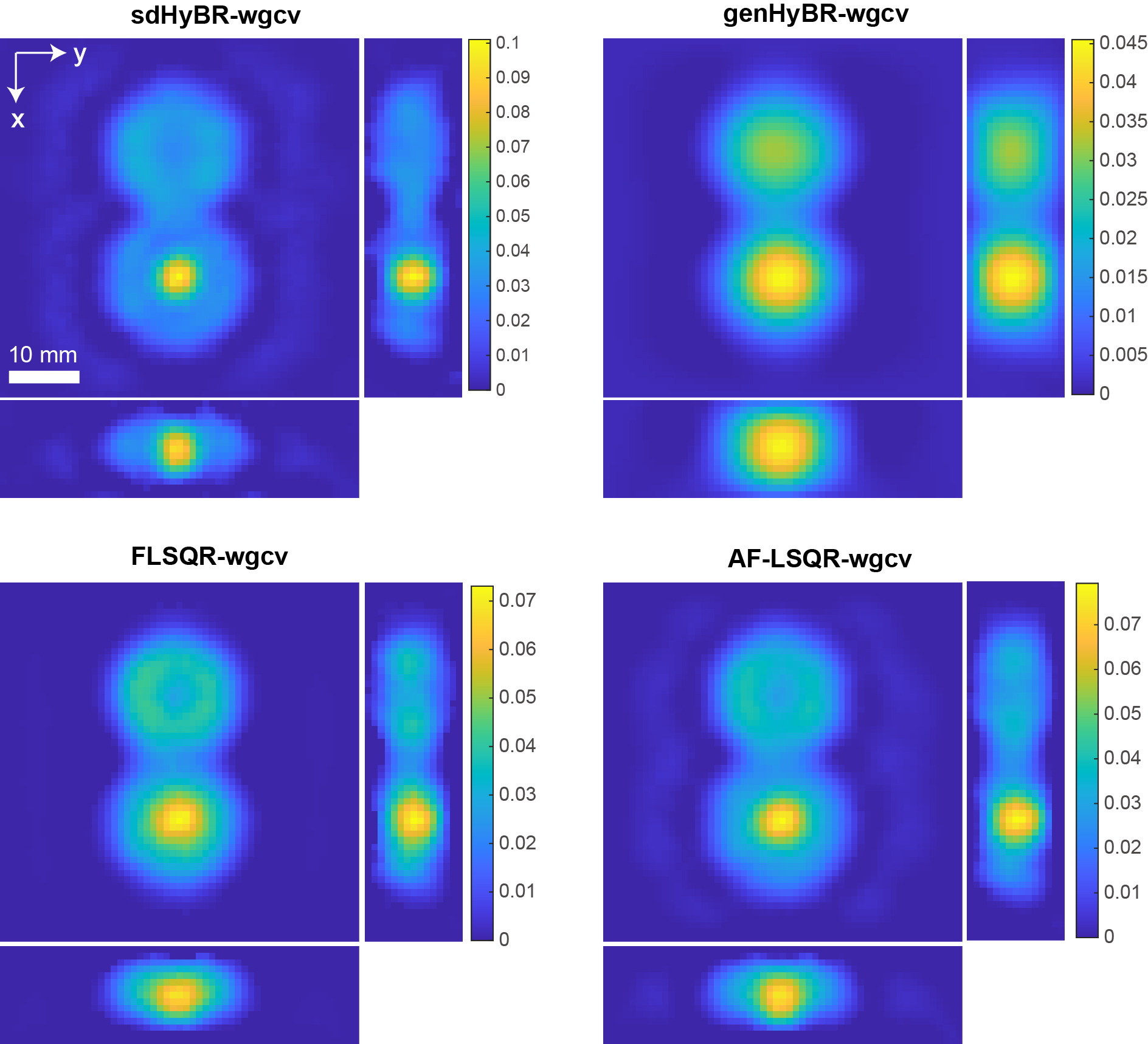}
   \centering
    \caption{FMT reconstruction using synthetic data: volumetric fluorescence distribution resulted from different reconstruction approaches. The three planes are visualized at the cross position of (38, 28, 8) mm. All results correspond to using WGCV-selected regularization parameters and are obtained at the stopping iterations.
   }
    \label{fig2_ex3}
\end{figure}

\begin{figure}[ht!]
   \includegraphics[width = 0.8\textwidth]{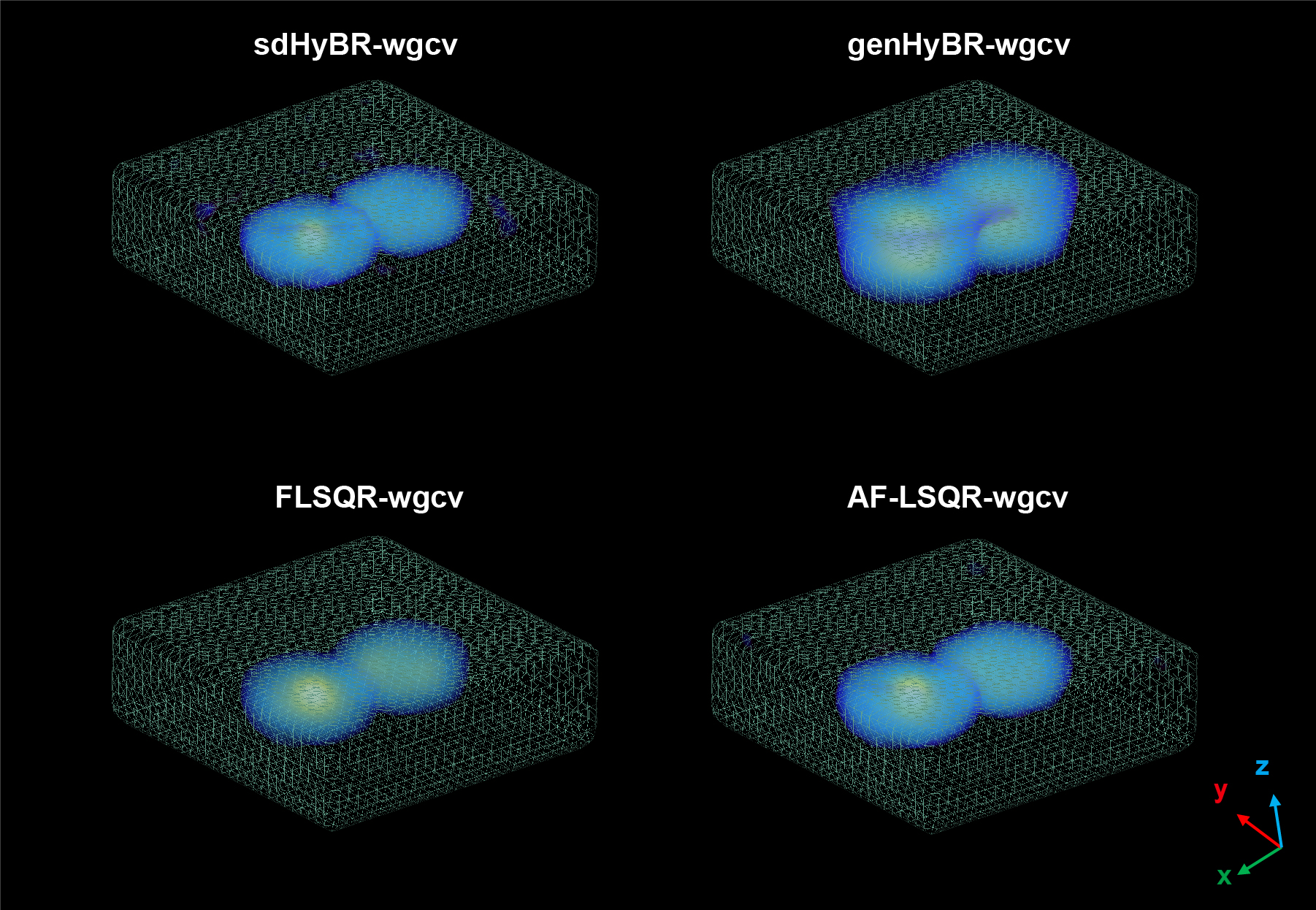}
   \centering
    \caption{3D visualization of the FMT reconstructions using synthetic data. All results correspond to using WGCV-selected regularization parameters and are obtained at the stopping iterations.}
    \label{fig3_ex3}
\end{figure}

\begin{figure}[ht!]
    \centering
    \begin{tabular}{ccc}
    \includegraphics[width = 0.48\textwidth]{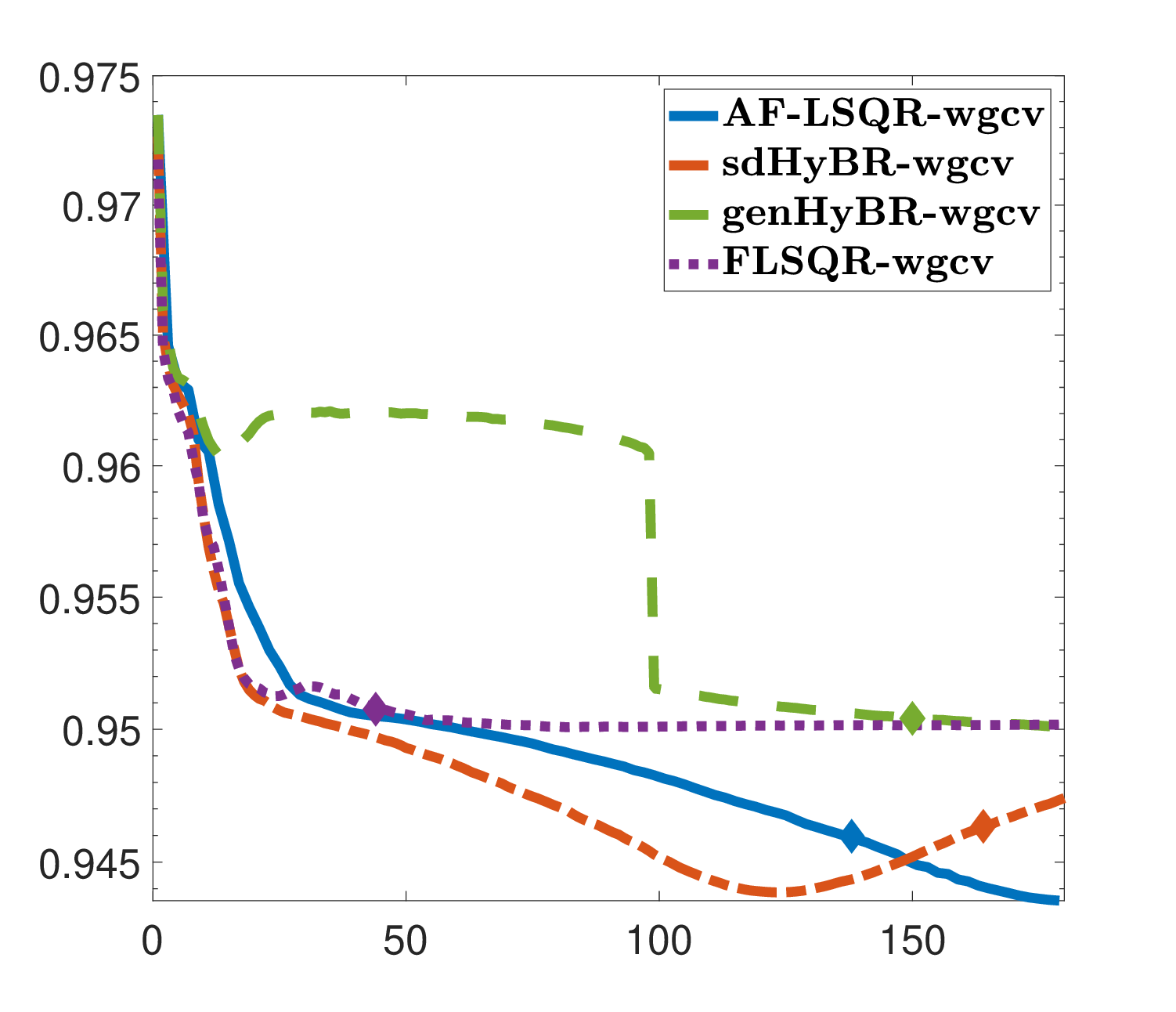} 
    &     \includegraphics[width = 0.5\textwidth]{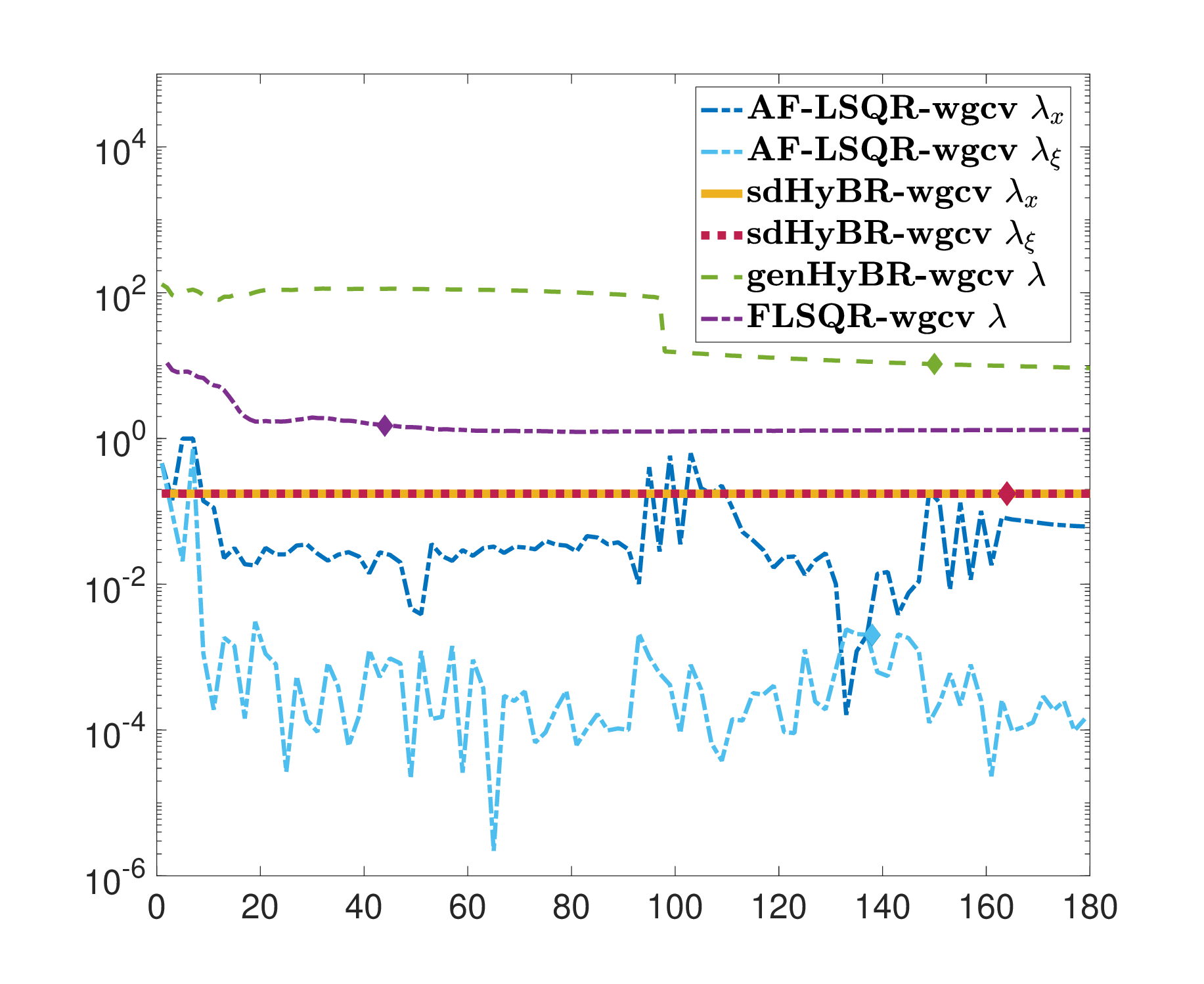} \\
         (a) Relative error norms & (b) Regularization parameters
    \end{tabular}
    \caption{(a) Relative error norms per iteration for the synthetic FMT problem in Figure \ref{fig1_ex3} for AF-LSQR, sdHyBR, genHyBR and FLSQR solvers. (b) Regularization parameters computed at each iteration using the WGCV-selected strategy. The flattening of the GCV function in \eqref{eq:GCVstop} is used to stop the iterations, and the stopping iteration is signaled in the plots with a marker. }
    \label{fig4_ex3}
\end{figure}

The fluorescence distributions obtained from all reconstruction methods are visualized in three cross-sections in Figure \ref{fig2_ex3}, and are consistently represented on a grid of $55\times55\times15$. 
The genHyBR method significantly amplifies the volume of the sparse fluorescence, leading to noticeable artifacts and unclear edges. The sdHyBR method introduces clearly visible ripple-like artifacts in smooth regions, which degrades the overall quality of its reconstruction. In comparison, both FLSQR and AF-LSQR exhibit a commendable balance between capturing sparse and smooth stochastic attributes, resulting in a more holistic reconstruction. However, compared with AF-LSQR, the reconstruction obtained using FLSQR still fails to recover appropriately the edges of the highly fluorescent region, posing challenges to practical applications. Overall, AF-LSQR stands out with the best qualitative performance, successfully capturing both sparse and smooth signals with superior image quality. Figure \ref{fig3_ex3} shows 3D visualizations of the compared reconstruction results, further demonstrating that the reconstruction of AF-LSQR is closer to the true fluorescence distribution compared to other approaches. Figure \ref{fig4_ex3} presents the relative reconstruction error norms and WGCV-selected regularization parameters, where the diamond-shaped markers denote the stopping iterations. Notably, AF-LSQR outperforms other algorithms in attaining the lowest relative error and ensuring stable convergence. While sdHyBR also achieves the minimal relative error, it exhibits semi-convergence. Moreover, one can observe in Figure \ref{fig4_ex3}(b) that the regularization parameters for sdHyBR barely change across the iterations, possibly indicating a bad performance of the regularization parameter choice criterion. Note that, since AF-LSQR and sdHyBR solve projections of different functional at each iteration, there is not a correlation between the regularization parameters that are appropriate for each method at each iteration.

\subsubsection{FMT reconstruction using real-world experimental data}
\begin{figure}[ht!]
   \centering
   \includegraphics[width = 1\textwidth]{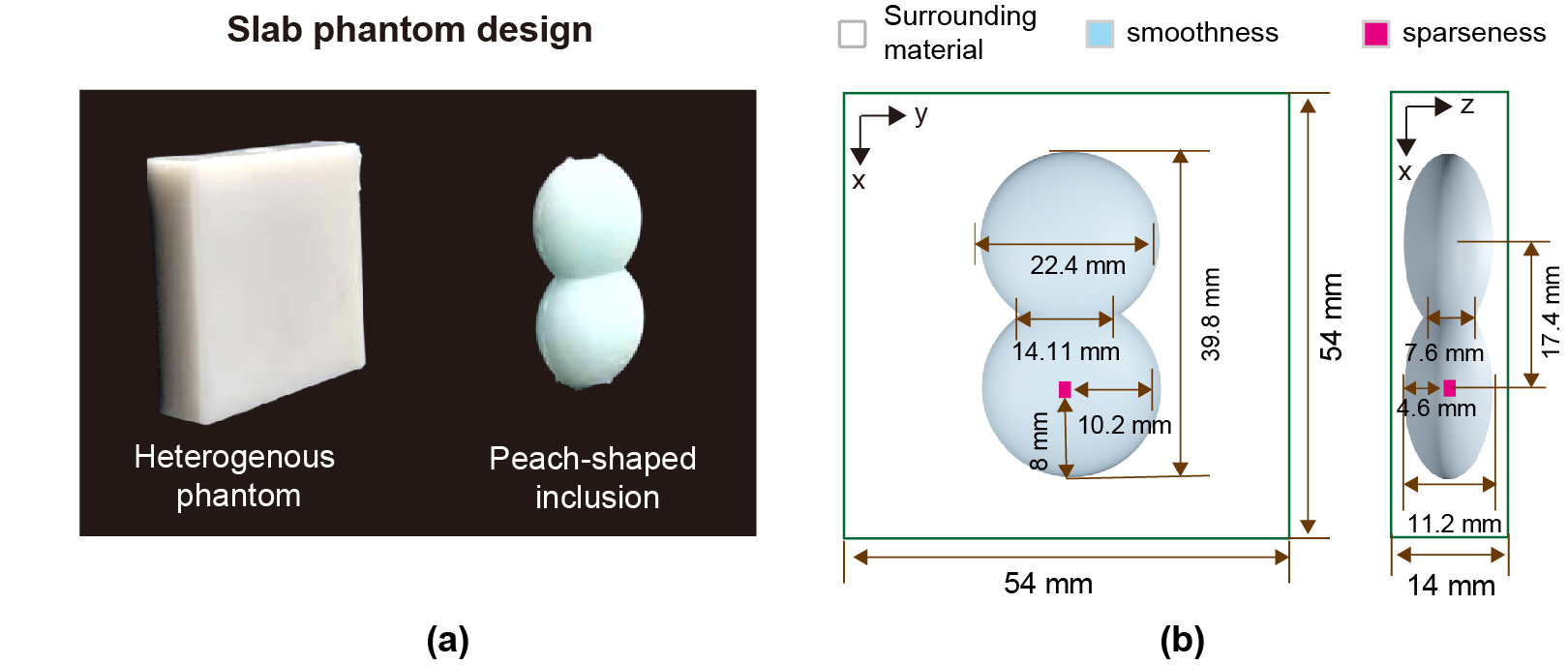}
    \caption{A real-world silicone phantom for a FMT experiment: (a) A silicone slab phantom used for evaluating our method. (b) The detailed design of the model. The optically heterogeneous phantom is designed by placing a large peanut-shaped region containing low-concentration fluorophore (blue) in the center. Through the center of the region, we insert a glass tube containing high concentration small volume fluorophore (red).}
    \label{fig1_ex4} 
\end{figure}

To test the potential use of AF-LSQR in real applications, we evaluate the performance of its reconstruction using a silicone phantom featuring similar optical properties to biological tissues. The design specifications, optical properties, and detailed imaging settings of the physical phantom precisely match those utilized in the simulated case outlined in Subsection \ref{subsec:fmt}. Specifically, to mimic the optical properties of biological tissues, the  fundamental structure of the phantom is composed of silicone (SYLGARD™ 184, DOW, CA, US), while the scattering and absorption coefficients modulated by introducing varying quantities of titanium dioxide (Ti$\mathrm{O}_\mathrm{2}$, Colins, Shanghai, China) particles and carbon black powder (Colins, Shanghai, China). The slab phantom is shown in Figure \ref{fig1_ex4}(a). Similar to the previous numerical case study, the fluorescence source consists of two components: 1) a large, peanut-shaped component containing low-concentration (0.0243 µmol ml$^{-1}$) fluorescence dye and 2) a small-volume saline-fluorophore mixture droplet with high concentration (1.9434 µmol ml$^{-1}$). The dimension, position and detailed design of the large-volume fluorescence inclusion are provided in Figure \ref{fig1_ex4}(b), whereas the fluorescence droplet is contained in a capillary glass tube (inner diameter: 1 mm, outer diameter: 1.80 mm) penetrating along the central axis of the large component. Data acquisition of the real phantom experiment was performed using a customized FMT system described in \cite{wu2023multifunctional}.

In this case study, we use the STIFT platform to construct the system matrix $A \in \mathbb{R}^{90197 \times 45375}$ and to pre-process the raw data accordingly to obtain the vector of measurements $b \in \mathbb{R}^{90197}$ used in the reconstructions. The construction of the covariance matrix $Q$ and the tolerance for the GCV function are identical to those of the simulated case described in Subsection \ref{subsec:fmt}. The maximum number of iterations was set to 40.

\begin{figure}[ht!]   
   \includegraphics[width = 0.8\textwidth]{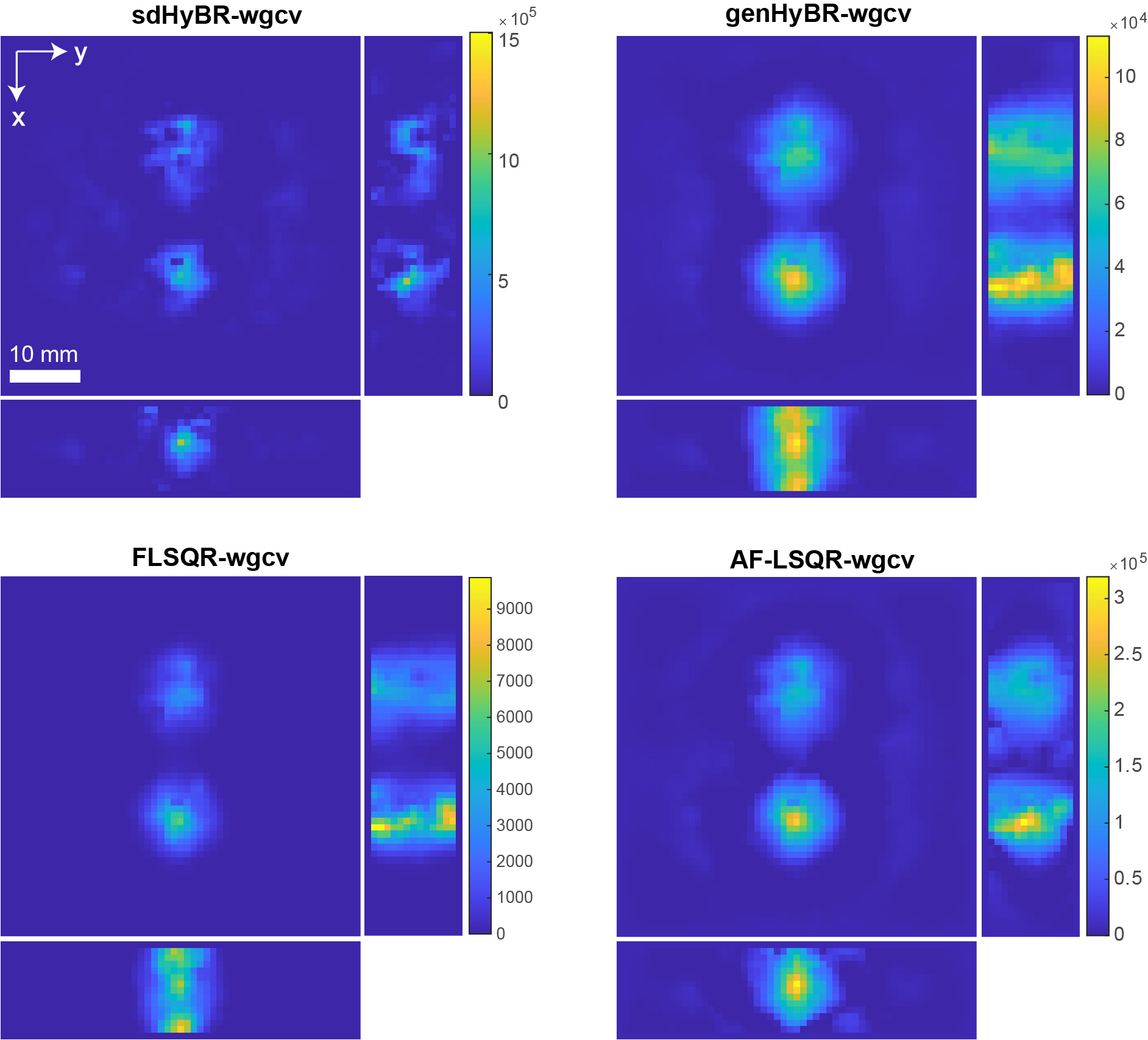}
   \centering
    \caption{FMT experimental study: 3D reconstructed fluorescence distribution resulted from different reconstruction approaches in the experiment study. The three planes are visualized at the cross position of (38, 28, 8) mm. All results correspond to using WGCV-selected regularization parameters and are obtained at the stopping iterations.}
    \label{fig2_ex4}
\end{figure}

\begin{figure}[ht!]
   \includegraphics[width = 0.8\textwidth]{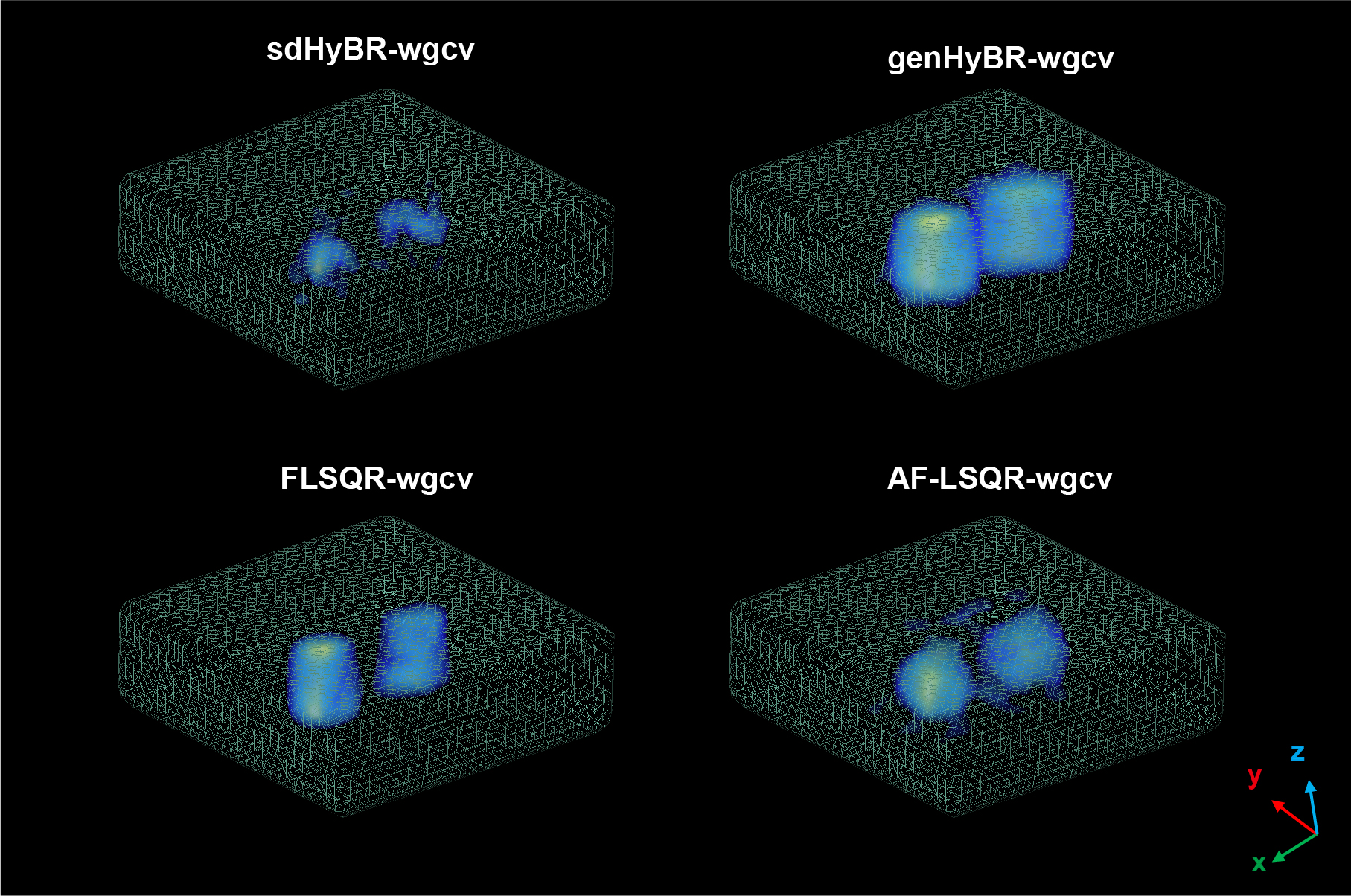}
   \centering
    \caption{FMT experimental study: Reconstruction results visualized in 3D for the FMT experimental study. All results correspond to using WGCV-selected regularization parameters and are obtained at the stopping iterations.}
    \label{fig3_ex4}
\end{figure}

Figure \ref{fig2_ex4} and Figure \ref{fig3_ex4} illustrate the reconstruction results from 2D cross-sections and 3D views respectively. One can observe that AF-LSQR successfully recovers an image displaying both smooth and sparse components. Although the boundary of the smooth component is perhaps overly smoothed, the location and shape of the phantom tumor can be clearly identified. In contrast, both sdHyBR and FLSQR are unable to accurately reconstruct the fluorescence patterns. For genHyBR the positional and edge information of the fluorescent regions is not as accurate as those of AF-LSQR, particularly for the high concentration volume. In Figure \ref{fig3_ex4}, we observe that the reconstructed morphology of the smooth regions from sdHyBR, genHyBR and FLSQR differs significantly from the true fluorescence distribution. The sdHyBR method completely fails to reconstruct the contour of the smooth region. In the reconstructions generated by genHyBR and FLSQR, the two lobes of the smooth region seem to approximate elliptical cylinders rather than the true peanut-like shape observed in \ref{fig1_ex4}(a). 

In summary, the results of this case study demonstrate that AF-LSQR can yield accurate results for complex, real-world FMT experiments. For this example, the reconstruction given by the proposed method provides an accurate positioning and clear edges for the sparse element without compromising the integrity of the smooth background, which is crucial for clinicians to identify conditions such as tumor locations.

\section{Discussion and conclusion}\label{Chapter:conclusions} 
In this paper we introduce two novel algorithms, namely AF-GMRES and AF-LSQR, designed to address large-scale Bayesian inverse problems where the solution contains two components with distinct stochastic characteristics (smoothness and sparsity). The new approaches are projection methods that (partially) solve a sequence of quadratic problems by projecting them into a (single) space of increasing dimensions, avoiding more classical but computationally expensive inner-outer schemes. The new methods are efficient and matrix-free, as they only require matrix-vector product with $A$ (and possibly $A^{\top}$) to compute the solution, making them suitable for large-scale problems. Moreover, AF-GMRES and AF-LSQR are naturally suited for the automatic selection of the different regularization parameters, and the projections are cleverly designed so that the parameters computed at each iteration for the appropriate projected problems have a correspondence to the regularization parameters of the full dimensional problem. Moreover, we provide theoretical results for the convergence of our methods when the regularization parameters are fixed. 

Extensive experiments, involving both numerical testing and a real-world problem, show that AF-GMRES and AF-LSQR are competitive with other hybrid methods. In particular, the results from the real FMT experiment demonstrate the accuracy, robustness, and potential of the new methods to solve real-world medical imaging problems. Future work includes extending and applying the proposed methods to more complex medical imaging problems such as using FMT to image the lymph nodes of mice.  Furthermore, AF-GMRES and AF-LSQR can be extended to solve problems including a variety of regularizers, including total variation and data-driven matrix-norm regularization. 

\section*{Acknowledgements and declarations of interest}
MSL gratefully acknowledges support from the Cambridge Mathematics of Information in Healthcare Hub (CMIH), University of Cambridge. This work was partially supported by the National Science Foundation program under grant DMS-2208294. Any opinions, findings, and conclusions or recommendations expressed in this material are those of the author(s) and do not necessarily reflect the views of the National Science Foundation. This work was also partially supported by National Natural Science Foundation of China under grant 12101406, grant 62105205, and Shanghai Science and Technology Innovation Program under grant 21YF1429100. \\
Declarations of interest: none.

\bibliographystyle{abbrv}
\bibliography{biblio}

\begin{thebibliography}{10}

\bibitem{Arridge_2014}
S.~R. Arridge, M.~M. Betcke, and L.~Harhanen.
\newblock Iterated preconditioned lsqr method for inverse problems on
  unstructured grids.
\newblock {\em Inverse Problems}, 30(7):075009, jun 2014.

\bibitem{arridge2009optical}
S.~R. Arridge and J.~C. Schotland.
\newblock Optical tomography: forward and inverse problems.
\newblock {\em Inverse problems}, 25(12):123010, 2009.

\bibitem{doi:10.1137/080716542}
A.~Beck and M.~Teboulle.
\newblock A fast iterative shrinkage-thresholding algorithm for linear inverse
  problems.
\newblock {\em SIAM Journal on Imaging Sciences}, 2(1):183--202, 2009.

\bibitem{JulianneSilvia}
J.~Chung and S.~Gazzola.
\newblock Flexible krylov methods for $\ell_p$ regularization.
\newblock {\em SIAM Journal on Scientific Computing}, 41(5):S149--S171, 2019.

\bibitem{hybrid_review}
J.~Chung and S.~Gazzola.
\newblock Computational methods for large-scale inverse problems: a survey on
  hybrid projection methods.
\newblock {\em arXiv preprint}, 2021.

\bibitem{sdecompose}
J.~Chung, J.~Jiang, S.~M. Miller, and A.~K. Saibaba.
\newblock Hybrid projection methods for solution decomposition in large-scale
  bayesian inverse problems.
\newblock 2022.

\bibitem{doi:10.1137/16M1081968}
J.~Chung and A.~K. Saibaba.
\newblock Generalized hybrid iterative methods for large-scale bayesian inverse
  problems.
\newblock {\em SIAM Journal on Scientific Computing}, 39(5):S24--S46, 2017.

\bibitem{https://doi.org/10.1002/cpa.20303}
I.~Daubechies, R.~DeVore, M.~Fornasier, and C.~S. Güntürk.
\newblock Iteratively reweighted least squares minimization for sparse
  recovery.
\newblock {\em Communications on Pure and Applied Mathematics}, 63(1):1--38,
  2010.

\bibitem{doi:10.1137/080731992}
J.~Demmel, L.~Grigori, M.~Hoemmen, and J.~Langou.
\newblock Communication-optimal parallel and sequential qr and lu
  factorizations.
\newblock {\em SIAM Journal on Scientific Computing}, 34(1):A206--A239, 2012.

\bibitem{gazzola2019ir}
S.~Gazzola, P.~C. Hansen, and J.~G. Nagy.
\newblock Ir tools: a matlab package of iterative regularization methods and
  large-scale test problems.
\newblock {\em Numerical Algorithms}, 81(3):773--811, 2019.

\bibitem{msl1}
S.~Gazzola, J.~G. Nagy, and M.~S. Landman.
\newblock Iteratively reweighted fgmres and flsqr for sparse reconstruction.
\newblock {\em SIAM Journal on Scientific Computing}, 43(5):S47--S69, 2021.

\bibitem{Krylov_props}
S.~Gazzola, P.~Novati, and M.~R. Russo.
\newblock On krylov projection methods and tikhonov regularization.
\newblock {\em Electronic Transactions on Numerical Analysis}, 44, 05 2014.

\bibitem{msl2}
S.~Gazzola and M.~Sabaté~Landman.
\newblock Krylov methods for inverse problems: Surveying classical, and
  introducing new, algorithmic approaches.
\newblock {\em GAMM-Mitteilungen}, 43(4):e202000017, 2020.

\bibitem{goldstein2009split}
T.~Goldstein and S.~Osher.
\newblock The split {B}regman method for l1-regularized problems.
\newblock {\em SIAM Journal on Imaging Sciences}, 2:323--343, 2009.

\bibitem{Hanke1993RegularizationMF}
M.~Hanke and P.~C. Hansen.
\newblock Regularization methods for large-scale problems.
\newblock {\em Surveys on Mathematics for Industry}, 3(4):253–315, 1993.

\bibitem{hilderbrand2010near}
S.~A. Hilderbrand and R.~Weissleder.
\newblock Near-infrared fluorescence: application to in vivo molecular imaging.
\newblock {\em Current opinion in chemical biology}, 14(1):71--79, 2010.

\bibitem{Huang2017}
G.~Huang, A.~Lanza, S.~Morigi, L.~Reichel, and F.~Sgallari.
\newblock Majorization-minimization generalized {K}rylov subspace methods for
  $l_p$-$l_q$ optimization applied to image restoration.
\newblock {\em BIT}, 57(2):351--378, Jun 2017.

\bibitem{lin2003near}
J.~Lin, C.~Gerbig, S.~Wofsy, A.~Andrews, B.~Daube, K.~Davis, and C.~Grainger.
\newblock A near-field tool for simulating the upstream influence of
  atmospheric observations: The stochastic time-inverted lagrangian transport
  (stilt) model.
\newblock {\em Journal of Geophysical Research: Atmospheres}, 108(D16), 2003.

\bibitem{liu2021data}
X.~Liu, A.~L. Weinbren, H.~Chang, J.~M. Tadi{\'c}, M.~E. Mountain, M.~E.
  Trudeau, A.~E. Andrews, Z.~Chen, and S.~M. Miller.
\newblock Data reduction for inverse modeling: an adaptive approach v1. 0.
\newblock {\em Geoscientific Model Development}, 14(7):4683--4696, 2021.

\bibitem{mieog2022fundamentals}
J.~S.~D. Mieog, F.~B. Achterberg, A.~Zlitni, M.~Hutteman, J.~Burggraaf, R.-J.
  Swijnenburg, S.~Gioux, and A.~L. Vahrmeijer.
\newblock Fundamentals and developments in fluorescence-guided cancer surgery.
\newblock {\em Nature Reviews Clinical Oncology}, 19(1):9--22, 2022.

\bibitem{miller2020geostatistical}
S.~M. Miller, A.~K. Saibaba, M.~E. Trudeau, M.~E. Mountain, and A.~E. Andrews.
\newblock Geostatistical inverse modeling with very large datasets: an example
  from the orbiting carbon observatory 2 (oco-2) satellite.
\newblock {\em Geoscientific Model Development}, 13(3):1771--1785, 2020.

\bibitem{nehrkorn2010coupled}
T.~Nehrkorn, J.~Eluszkiewicz, S.~C. Wofsy, J.~C. Lin, C.~Gerbig, M.~Longo, and
  S.~Freitas.
\newblock Coupled weather research and forecasting--stochastic time-inverted
  lagrangian transport (wrf--stilt) model.
\newblock {\em Meteorology and Atmospheric Physics}, 107:51--64, 2010.

\bibitem{ntziachristos2006fluorescence}
V.~Ntziachristos.
\newblock Fluorescence molecular imaging.
\newblock {\em Annu. Rev. Biomed. Eng.}, 8:1--33, 2006.

\bibitem{10.7551/mitpress/3206.001.0001}
C.~E. Rasmussen and C.~K.~I. Williams.
\newblock {\em {Gaussian Processes for Machine Learning}}.
\newblock The MIT Press, 11 2005.

\bibitem{ren2019smart}
W.~Ren, H.~Isler, M.~Wolf, J.~Ripoll, and M.~Rudin.
\newblock Smart toolkit for fluorescence tomography: simulation,
  reconstruction, and validation.
\newblock {\em IEEE Transactions on Biomedical Engineering}, 67(1):16--26,
  2019.

\bibitem{doi:10.1137/15M1037925}
R.~A. Renaut, S.~Vatankhah, and V.~E. Ardestani.
\newblock Hybrid and iteratively reweighted regularization by unbiased
  predictive risk and weighted gcv for projected systems.
\newblock {\em SIAM Journal on Scientific Computing}, 39(2):B221--B243, 2017.

\bibitem{renaut2017hybrid}
R.~A. Renaut, S.~Vatankhah, and V.~E. Ardestani.
\newblock Hybrid and iteratively reweighted regularization by unbiased
  predictive risk and weighted gcv for projected systems.
\newblock {\em SIAM Journal on Scientific Computing}, 39(2):B221--B243, 2017.

\bibitem{IRN2}
B.~Wohlberg and P.~Rodriguez.
\newblock An efficient algorithm for sparse representations with $\ell^p$ data
  fidelity term.
\newblock {\em In Proceedings of the 4th IEEE Andean Technical Conference
  (AN-DESCON)}, 2008.

\bibitem{3883}
S.~Wright, R.~Nowak, and M.~A.~T. Figueiredo.
\newblock Sparse reconstruction by separable approximation.
\newblock {\em IEEE Transactions on Signal Processing}, 57(7):--, July 2009.

\bibitem{wu2023multifunctional}
Y.~Wu, S.~Gao, L.~Li, J.~Zhang, Q.~Hu, X.~Lou, X.~Zhu, J.~Jiang, and W.~Ren.
\newblock Multifunctional optical tomography system combining surface
  extraction and 3d fluorescence reconstruction.
\newblock In {\em Second Conference on Biomedical Photonics and Cross-Fusion
  (BPC 2023)}, volume 12753, pages 108--114. SPIE, 2023.

\end{thebibliography}

\end{document}